\documentclass[a4paper,draft]{amsproc}
\usepackage{amssymb}
\usepackage{amscd}
\theoremstyle{plain}
 \newtheorem{thm}{Theorem}[section]
 \newtheorem{prop}{Proposition}[section]
 \newtheorem{lem}{Lemma}[section]
 \newtheorem{cor}{Corollary}[section]
\theoremstyle{definition}
 \newtheorem{exm}{Example}[section]
 \newtheorem{rem}{Remark}[section]
\newtheorem{dfn}{Definition}[section]

\numberwithin{equation}{section}

\renewcommand{\leq}{\leqslant} \renewcommand{\geq}{\geqslant}

\newcommand{\R}{\mathbb{R}}

\newcommand{\ddim}{\mathrm{ddim\;}}
\newcommand{\dind}{\mathrm{dind\;}}
\newcommand{\corank}{\mathrm{corank\;}}
\newcommand{\tr}{\mathrm{tr}}

\newcommand{\T}{\mathbb{T}}
\newcommand{\g}{\mathfrak{g}}

\newcommand{\ad}{\mathrm{ad}}
\DeclareMathOperator{\Ad}{\mathrm{Ad}}

\DeclareMathOperator{\diag}{\mathrm{diag}}
\DeclareMathOperator{\rank}{\mathrm{rank}}
\DeclareMathOperator{\pr}{\mathrm{pr}}

\newcommand{\D}{\mathfrak{ d}}

\title[SYSTEMS OF HESS--APPEL'ROT TYPE AND ZHUKOVSKII PROPERTY]{SYSTEMS OF HESS--APPEL'ROT TYPE AND ZHUKOVSKII PROPERTY}

\author[Dragovi\'c, Gaji\'c, Jovanovi\'c]{\bfseries Vladimir Dragovi\'c$^{1,2}$,
Borislav Gaji\'c$^1$, Bo\v zidar Jovanovi\'c$^1$}

\address{$^1$\,Mathematical Institute \\Serbian Academy of Sciences and Arts \\Kneza Mihaila 36, 11000 Belgrade\\Serbia}
\address{$^2$\, GFM, University of Lisbon, Portugal}
\address{\rm e-mail addresses: vladad@mi.sanu.ac.rs, gajab@mi.sanu.ac.rs,
bozaj@mi.sanu.ac.rs}

\begin{document}

\vspace*{18mm} \setcounter{page}{1} \thispagestyle{empty}

\maketitle

 \centerline{\it Dedicated to the memory
of Professor Novica Bla\v zi\' c (1959-2005)}

\vspace{15mm}

\begin{abstract}
We start with a review of a class of systems with invariant
relations, so called {\it systems of Hess--Appel'rot type} that
generalizes the classical Hess--Appel'rot rigid body case.  The
systems of Hess--Appel'rot type have remarkable property: there
exists a pair of compatible Poisson structures, such that a system
is certain Hamiltonian perturbation of an integrable
bi-Hamiltonian system. The invariant relations are Casimir
functions of the second structure. The systems of Hess--Appel'rot
type carry an interesting combination of both integrable and
non-integrable properties.

Further, following integrable line, we study partial reductions
and systems having what we call the {\it Zhukovskii property}:
These are Hamiltonian systems on a symplectic manifold $M$ with
actions of two groups $G$ and $K$; the systems are assumed to be
$K$--invariant and to have invariant relation $\Phi=0$ given by
the momentum mapping of the $G$--action, admitting two type of
reductions, a reduction to the Poisson manifold $P=M/K$ and a
partial reduction to the symplectic manifold $N_0=\Phi^{-1}(0)/G$;
final and crucial assumption is that the partially reduced system
to $N_0$ is completely integrable. We prove that the Zhukovskii
property is a quite general characteristic of  systems of
Hess--Appel'rote type.  The partial reduction neglects the most
interesting and challenging part of the dynamics of the systems of
Hess--Appel'rot type - the non-integrable part, some analysis of
which  may be seen as a reconstruction problem.

We show that an integrable system, the magnetic pendulum on the
oriented Grassmannian $Gr^+(n,2)$ has a natural interpretation
within Zhukovskii property and that it is equivalent to a partial
reduction of certain system of Hess--Appel'rot type. We perform a
classical and algebro-geometric integration of the system in
dimension four, as an example of a known isoholomorphic system -
the Lagrange bitop.

The paper presents a lot of examples of systems of Hess--Appel'rot
type, giving an additional argument in favor of further study of
this class of systems.
\end{abstract}

\newpage
\tableofcontents

\section{Introduction}

Historically, the Hess-Appel'rot system as a classical rigid-body
system, appeared (see \cite{He}) just a year after the celebrated
Kowalevski  1889 paper \cite{Kow}, and its immediate popularity
had been connected with its relationship with the Kowalevski
paper. Kowalevski  started, as we know, from a careful analysis of
the solutions of the Euler and the Lagrange case of rigid-body
motion and formulated a problem {\it of describing the parameters
$(A, B, C, x_0, y_0, z_0)$, for which the Euler--Poisson equations
have a general solution in a form of uniform functions only with
moving poles as singularities.} Here, $I=\diag (A, B, C)$
represents the inertia operator, and $\chi=(x_0, y_0, z_0)$ is the
centre of mass of the rigid body.

Then, in \S 1 of \cite{Kow}, after some necessary  conditions were
formulated, she descovered a new case, which is now known as the
Kowalevski case. The last  case  was, according to Kowalevski, a
unique possible beside the cases of Euler and Lagrange. However,
considering the situation where all momenta of inertia are
different, Kowalevski came to a relation analogue to the following
(see \cite{Go}):
$$
x_0\sqrt{A(B-C)}+y_0\sqrt{B(C-A)}+z_0\sqrt{C(A-B)}=0,
$$
and she concluded that $x_0=y_0=z_0$, which represented  the Euler
case.

However, Appel'rot  noticed few years later, that the last relation
admitted one more case, not observed in \cite{Kow}:
$$
x_0\sqrt{A(B-C)}+z_0\sqrt{C(A-B)}=0,\quad y_0=0,
$$
where he  assumed $A>B>C$.  Such intriguing position corresponding
to the possible mistake  in the Kowalevski paper, made the
Hess-Appel'rot systems very attractive for leading Russian
mathematicians from the end of XIX century as a possible
counterexample. But, after a few years, Nekrasov and Lyapunov proved
that the Hess-Appel'rot systems didn't satisfy the condition
investigated by Kowalevski, which means that conclusion of \S 1 of
\cite{Kow} {\it was correct}.

And, from that moment until very recently, the Hess-Appel'rot
systems were basically left aside, even in modern times, when new
methods of inverse problems, Lax representations,  finite-zone
integrations were applied to almost all known classical systems.

A modern theory of systems of Hess-Appel'rot type has been
developed in \cite {DrGa, DrGa1, DrGa3, DrGa2}. It started with a
construction of a Lax representation for the Hess-Appel'rot system
in \cite{DrGa}, see the Proposition \ref{proplapar} below.
Generalization of this Lax pair in four-dimensional case led to
construction of a new integrable rigid-body system in \cite{DrGa},
called {\it the Lagrange bitop}. Algebro-geometric integration
procedure of the Lagrange bitop has been performed in
\cite{DrGa1}. It brought to the discovery of a new class of
integrable systems, which was named {\it isoholonomic systems} in
\cite{DrGa1}. Higher dimensional generalizations of the classical
rigid-body Hess-Appel'rot systems have been constructed in
\cite{DrGa3, DrGa2} as certain perturbations of the isoholonomic
integrable systems. Finally, in \cite{DrGa2} after detailed
analysis of infinite set of new examples in arbitrary dimensions,
the theory of systems of Hess-Appel'rot type has been settled down
providing an axiomatic, general and abstract approach, see also
the Section \ref{sec:axiom}.

According to this theory, the systems of Hess-Appel'rot type form a
class of dynamical systems, obtained as certain perturbations of
integrable, bi-Hamiltonian systems which carry an interesting
combination of both integrable and non-integrable properties.

Suppose a bi-Poisson structure
$\{\cdot,\cdot\}_1+\lambda\{\cdot,\cdot\}_2$ is given, with an
integrable, bihamiltonian system with the Hamiltonian $H_0$
corresponding to the first structure. Further, let $f_1,\dots, f_k$
be the commuting integrals of the system $(H_0,\{\cdot,\cdot\}_1)$,
which are Casimirs for the second structure $\{\cdot,\cdot\}_2$.
Then, a system of Hess-Appel'rot type is Hamiltonian with respect to
the first structure with a Hamiltonian
$$
H=H_0+\sum_{l=1}^k J_lb_lf_l,
$$
where $J_l$ are constants and $b_l$ are certain functions on the
phase space. The invariant relations are
$$
f_l=0,\quad l=1,\dots, k.
$$
Thus, the invariant manifolds are symplectic leaves of the second
Poisson structure.

As perturbations, they are global and not just small perturbations
as it is usually the case in the study of non-integrable
perturbations of integrable systems. A balance between integrable
and nonintegrable properties is obtained by the choice of
perturbations. The system of Hess-Appel'rot type is Hamiltonian with
respect to the first Poisson structure, but perturbations and
invariant relations are defined by Casimirs of the second structure.

For the classical Hess-Appel'rot case, its integrable part made it
close to the Ko\-wa\-le\-vski study while its nonintegrable side
finally disqualified it as a possible counterexample for the
Kowalevski statement. Classical integration of its integrable part
one may find also in the book of Golubev \cite{Go}, while
algebro-geometric integration has been performed in \cite{DrGa}.
It was Zhukovskii who observed (see \cite{Zh}), that after certain
reduction, the classical Hess-Appel'rot system reduces to the
completely integrable system of spherical pendulum (see the
Subsection \ref{ha-sp} below). This observation of Zhukovskii
motivated us to introduce the notion of {\it Zhukovskii property},
see the Section see Section \ref{sec:reduction}.

For the four-dimensional generalization of the Hess-Appel'rot
system, detailed separation of integrable and nonintegrable part
has been done in two ways in \cite{DrGa3, DrGa2}, both classically
and algebro-geometrically. Moreover, two integration procedures
for the integrable part have been performed in all details in
\cite{DrGa2}.

From these integrations, one can see that a completely integrable
system in the same space requires one integration more in
three-dimensional case and in the four-dimensional case it would
require two integrations more.

Reviews of these results can be found in \cite{Drag} and \cite{Gaj}.

Following \cite{DrGa2}, the interest for the Zhukovskii property
in the modern context has been expressed in \cite{Jo2, Jo3}, where
a study of partial reductions has been developed further, see
Section \ref{sec:reduction}.

One of the aims of the present paper is to provide a systematic
study of the Zhukovski property from a general point of view of
the axioms of the systems of Hess-Appel'rot type. It appears that
the Zhukovskii property is a quite general characteristic of the
systems of Hess-Appel'rot type and together with the partial
reduction it traces well the integrable part of a system of
Hess-Appel'rot type, see Sections \ref{sec:reduction},
\ref{sec:partial}. In the same time, it neglects totally the most
interesting and challenging part of the dynamics of the systems of
Hess-Appel'rot type - the non-integrable part. This blindness to
the non-integrable part, makes the partial reduction being
important but of a limited range and domain in complete
understanding and studying of the systems of Hess-Appel'rot type.

Nevertheless, a completely  integrable system, the magnetic
pendulum on $Gr^+(4,2)$ which has been introduced and studied in
\cite{BJ5} within the study of magnetic flows on homogeneous
spaces, appears to have an interpretation within Zhukovski
property: we show that it can be obtained as a partial reduction
of a certain system of Hess-Appel'rot type. Moreover, it appears
that the magnetic pendulum is equivalent to a very simple instant
of the Lagrange bitop from \cite{DrGa1}. Thus, the magnetic
pendulum is an example of an isoholonomic system. So, the
integration techniques of \cite{DrGa1} and \cite{DrGa2} can bee
applied directly to the magnetic pendulum. Complete integration of
the magnetic pendulum in classical and algebro-geometric manner is
the second main aim of this paper and it is performed in the
Section \ref{sec:pendulum}.

The paper provides a lot of examples of systems of Hess-Appel'rot
type, motivated from \cite{DrGa3, DrGa2, Jo2, Jo3, DGJ} and from
references therein, see the Sections \ref{sec:examples1} and
\ref{sec:examples2}. Such a rich set of examples is additional
argument in favor of further study of the systems of
Hess-Appel'rot type.

At the end of the paper, we collected for a reader's sake all
necessary notions from the theory of Hamiltonian systems and their
reductions, see the Appendix.

\section{Classical three-dimensional Hess-Appel'rot system}
\label{sec:classicalHA}

The Euler-Poisson equations of the motion of a heavy rigid body in
the moving frame are \cite{Go}:
\begin{equation}\label{eq:EPeq}
\aligned \dot {\bold M}&={\bold M}\times {\bold \Omega}+
{\bold \Gamma}\times{\bold \chi},\\
\dot {\bold \Gamma}&={\bold \Gamma}\times{\bold \Omega}\\
{\bold \Omega}&=\tilde{J}{\bold M},\ \ \tilde{J}=\diag(\tilde{J}_1,
\tilde{J}_2, \tilde{J}_3),
\endaligned
\end{equation}
where $\bold M$ is the kinetic momentum vector, $\bold \Omega$ the
angular velocity, $\tilde{J}$ a diagonal matrix, the inverse of
inertia operator, $\bold \Gamma$ a unit vector fixed in the space
and $\bold \chi$ is the radius vector of the centre of masses.

It is well known (see for example \cite{Go}) that equations
(\ref{eq:EPeq}) have three integrals of motion:
\begin{equation}\label{eq:intrigbody}
F_1=\frac 12\langle {\bold M},
{\bold\Omega}\rangle+\langle{\bold\Gamma},
{\bold\chi}\rangle,\qquad F_2=\langle {\bold M},
{\bold\Gamma}\rangle,\qquad F_3=\langle {\bold \Gamma},
{\bold\Gamma}\rangle=1.
\end{equation}

Thus, for complete integrability, one integral more is necessary
\cite{Go}. Let $\tilde{J}_1<\tilde{J}_2<\tilde{J}_3$ and ${\bold
\chi}=(x_0, y_0, z_0)$. Hess in \cite{He} and Appel'rot in
\cite{App} found that if the inertia momenta and the radius vector
of the centre of masses satisfy the conditions
\begin{equation}\label{eq:HAdef}
\aligned
y_0&=0\\
x_0\sqrt{\tilde{J}_3-\tilde{J}_2}&+z_0\sqrt{\tilde{J}_2-\tilde{J}_1}=0,
\endaligned
\end{equation}
then the surface
\begin{equation}\label{eq:HAinv}
F_4=M_1x_0+M_3z_0=0
\end{equation}
is invariant.

The compact connected components of the regular invariant sets
given by \eqref{eq:intrigbody}, \eqref{eq:HAinv} are tori, but not
with quasi-periodic dynamics. The classical and algebro-geometric
integration can be found in \cite{Go} and \cite{DrGa},
respectively. It is shown that the equations of the motion reduce to one
 elliptic integral and one Riccati differential equation.

{\it The Zhukovskii  geometric interpretation of the conditions}
\eqref{eq:HAdef} \cite{Zh, Leim}. Let us consider the ellipsoid
$$
\frac {M_1^2}{\tilde{J}_1}+\frac {M_2^2}{\tilde{J}_2}+\frac
{M_3^2}{\tilde{J}_3}=1,
$$
and the plane containing the middle axis and  intersecting the
ellipsoid at a circle. Denote by $l$ the normal to the plane,
which passes through the fixed point $O$. Then the condition
\eqref{eq:HAdef} means that the centre of masses  lies on the line
$l$.

Having this interpretation in mind, we choose a basis of moving
frame  such that the third axis is $l$, the second one is directed
along the middle axis of the ellipsoid, and the first one is chosen
according to the orientation of the orthogonal frame. In this basis
(see \cite{BoMa}), the particular integral (\ref{eq:HAinv}) becomes
$$
F_4=M_3=0,
$$
and the matrix $\tilde{J}$ and mass centre $\chi$ obtain the form:
\begin{equation}\label{Zhukovski}
J=\left(\begin{array}{ccc} J_1&0&J_{13}\\
                0&J_1&0\\
                J_{13}&0&J_3
        \end{array}\right), \qquad {\bold \chi}=(0, 0, z_0).
\end{equation}
This will serve us as a motivation for a definition of the
four-dimensional Hess-Appel'rot system.

A three-dimensional Lagrange top is defined by the Hamiltonian:
$$
H_L=\frac 12
\left(\frac{M_1^2+M_2^2}{I_1}+\frac{M_3^2}{I_3}\right)+z_0\Gamma_3,
$$
according to the standard Poisson structure
$$
\{M_i, M_j\}_1=-\epsilon_{ijk}M_k,\,\,
\{M_i,\Gamma_j\}_1=-\epsilon_{ijk}\Gamma_k,\,\,\{\Gamma_i,
\Gamma_j\}=0
$$
on the Lie algebra $e(3)$. It is also well-known that
three-dimensional Lagrange top is Hamiltonian in another Poisson
structure, compatible with first one. This structure is defined by:
$$
\{\Gamma_i,\Gamma_j\}_2=-\epsilon_{ijk}\Gamma_k,\,\,\{M_1,
M_2\}_2=1,
$$
and the corresponding Hamiltonian is:
$$
\tilde H_L=(a-1)M_3\left(\frac12
(M_1^2+M_2^2)+\Gamma_3\right)+M_1\Gamma_1+M_2\Gamma_2+M_3\Gamma_3
$$
where $I_1=1, I_3=a, z_0=1$. Casimir functions in the second
structure are $\Gamma_1^2+\Gamma_2^2+\Gamma_3^2$ and $M_3$.

Let us observe that the Hamiltonian for the three-dimensional
Hess-Appel'rot case is a quadratic deformation of Hamiltonian
$H_L$ of the Lagrange top:
$$
H_{HA}=H_L+J_{13}M_1M_3.
$$
The function $M_3$, which gives the invariant relation for the
Hess-Appel'rot case, is a Casimir function of the second Poisson
structure.

Let us mention that the Lax representation for Hess-Appel'rot system
\eqref{eq:EPeq}, \eqref{eq:HAdef}, \eqref{eq:HAinv} is constructed
in \cite{DrGa}:
\begin{prop}{\label{proplapar}}{\rm (\cite{DrGa})}
 On invariant manifold given by the invariant relation \eqref{eq:HAinv},
the equations of Hess-Appel'rot system are equivalent to the
matrix equation
$$
\dot{L}(\lambda)=[L(\lambda), A(\lambda)]
$$
where $ L(\lambda)=\lambda^2C+\lambda M+\Gamma$,
$A(\lambda)=\lambda\chi+\Omega. $ \end{prop} Here we denoted with
$M\in so(3)$ antisymmetric matrix that corresponds to the vector
$\bold{M}\in\mathbb{R}^3$ due to correspondence
$M_{ij}=-\epsilon_{ijk}M_{k}$ (and similar for $C, \Gamma, \Omega,
\chi$), and $C=\frac{1}{\tilde{J}_2}\chi$. Another Lax
representation for the Hess-Appel'rot system is given in \cite{DrGa4}. Using
it a sort of separation of variables for Hess-Appel'rot system is discussed there. By putting
$n=3$ in  \eqref{LA_pair} we get yet another Lax representation (see \cite{Jo2}).

\section{General systems of Hess-Appel'rot type. Axiomatic
approach}\label{sec:axiom}

\subsection{The first set of axioms: general Poisson settings}

Suppose a Poisson manifold $(M^{2n}, \{\cdot,\cdot\})$ is given,
together with $k+1$ functions $H, f_1,\dots, f_k \in C^{\infty}(M)$,
such that

\begin{itemize}
\item[(A1)]
$$
\{H, f_i\}= \sum_{j=1}^k a_{ij}f_j, \quad a_{ij}\in
C^{\infty}(M),\quad i, j=1,\dots ,k;
$$
\item[(A2)]
$$
\{f_i,f_j\}=0,\quad i,j=1,\dots ,k.
$$

A more general case can be obtained by replacing condition (A2)
with \item[(A2')]
$$
\{f_i,f_j\}=\sum_{l=1}^k d_{ij}^l f_l,\quad d_{ij}^l=const,\ \
i,j=1,...,k.
$$

In this case, the algebra of invariant relations is a noncommutative
Lie algebra.

Starting from the Hamiltonian system $(M, H_0)$ with $k$ integrals
in involution $f_1,\dots, f_k$, choosing functions $b_j\in
C^{\infty}(M), \quad j=1,\dots ,k$, one comes to a restrictively
integrable system:

\item[(HP)] {\bf (Hamiltonian perturbation axiom)}

The system $(M, H)$ where
$$
H=H_0 + \sum _{j=1}^kb_j f_j,
$$
will be called a Hamiltonian perturbation. It satisfies (A1) with
$$
a_{ij}=\{b_j, f_i\}, \quad i,j=1,\dots, k.
$$

\item[(BP)] {\bf (Bi-Poisson axiom)}
 {\it There exist a pair of compatible Poisson
structures, such that the system is Hamiltonian with respect to the
first structure, having the Hamiltonian of the form {\rm (HP)}, such
that $f_i$ are Casimir functions with respect to the second
structure.}
\end{itemize}

The invariant relations define symplectic leaves with respect to
the second structure, and the system is Hamiltonian with respect
the first one.

\subsection{The second set of axioms: Kowalevski property}
To get the right choice of axioms, we have to turn back to the
Kowalevski analysis. First, we are going to introduce some general
notions, see \cite{Koz}.

Suppose a system of ODEs of the form
\begin{equation}\label{eq:ODE}
\dot z_i=f_i(z_1,\dots,z_n), \quad i=1,\dots, n,
\end{equation}
is given and  there exist positive integers $g_i, \quad i=1,\dots,
n$, such that
$$
f_i(a^{g_1}z_1,\dots,a^{g_n}z_n)=a^{g_i+1}f_i(z_1,\dots, z_n), \quad
i=1,\dots, n.
$$
Then the system (\ref{eq:ODE}) is {\it quasi-homogeneous} and
numbers $g_i$ are {\it exponents of quasi-homogeneity}. Then, for
any complex solution $C=(C_1,\dots,C_n)$ of the system of
algebraic equations:
\begin{equation}\label{eq:Kowmat}
-g_iC_i=f_i(C_1,\dots,C_n), \quad i=1,\dots, n,
\end{equation}
one can define {\it the Kowalevski matrix} $K=K(C)=[K^i_j(C)]$:
$$
K^i_j(C)=\frac{\partial f_i}{\partial z_j}(C)+g_i\delta^i_j.
$$

Eigen-values of the Kowalevski matrix are called {\it the
Kowalevski exponents}. This terminology was introduced in
\cite{Yosh}. In last twenty years, heuristic and theoretical
methods in application of Kowalevski matrix and Kowalevski
exponents in study of integrability and nonintegrability have been
actively developing, see for example \cite{AvM2, AvMV, Koz}. But,
the notion of Kowalevski matrix and Kowalevski exponents were
introduced by Kowalevski herself in \cite{Kow}. The criterion she
used (see \cite{Kow}, p. 183, l. 15-22) to detect a system which
is now known as the Kowalevski top, can be formulated in Yoshida
terminology as:

{\bf Kowalevski condition (Kc).} {\it The $6\times 6$ Kowalevski
matrix should have} five different positive integer {\it Kowalevski
exponents}.

Now we return to the study of systems of Hess-Appel'rot type. The
systems we have constructed {\it are quasi-homogeneous}. Exponents
of each $M$ variable are $g=1$, and for any $\Gamma$ they are equal
to two. We are going now to calculate Kowalevski exponents for the
Hess-Appel'rot systems.

{\bf Three-dimensional Hess-Appel'rot case.} Let us denote $(M_1,
M_2, M_3, \Gamma_1, \Gamma_2$, $ \Gamma_3)$ by $(z_1,\dots, z_6)$.
Then the Euler-Poisson equations take the form \eqref{eq:ODE} with
$$
\aligned
f_1&=(J_3-J_1)z_2z_3+J_{13}z_1z_2+z_5;\\
f_2&=(J_3-J_1)z_1z_3+J_{13}(z_3^2-z_1^2)-z_4;\\
f_3&=J_{13}z_2z_3;\\
f_4&=J_3z_5z_3-J_1z_2z_6+J_{13}z_1z_5;\\
f_5&=-J_3z_3z_4+J_1z_1z_6+J_{13}(z_3z_6-z_1z_4);\\
f_6&=J_3z_2z_4-J_1z_1z_5-J_{13}z_3z_5;\\
\endaligned
$$
and $g_i=1, \quad i=1, 2, 3$ and $g_i= 2,\quad i= 4, 5, 6.$ The
invariant relation corresponds to the constraint $c_3=0$. So, we
are looking for solutions $(c_1,c_2,0,c_4,c_5,c_6)$ of the system
of the form \eqref{eq:Kowmat}. One can easily get
$c_4=-J_{13}c_1^2+c_2$, $c_5=-c_1(1+J_{13}c_2)$,
$c_6=-(c_1^2+c_2^2)/2$. Then, for $c_1\ne 0$, we get four possible
solutions for $(c_1, c_2)$ divided into two pairs: $(\pm
i/J_{13},-1/J_{13})$ and $(\pm 2i/J_{13},-2/J_{13})$. The
Kowalevski exponents are
$$
(-1, -2, 2, 4, 3, 3),\qquad (-1, 1, 3, 2, 2, 2),
$$
respectively.

Thus, it can easily be seen that classical Hess-Appel'rot system
doesn't satisfy exactly the Kowalevski condition (Kc), although it
is quite close to.

Thus, using into account properties of Kowalevski exponents of
algebraically-inte\-gra\-ble Hamiltonian systems, we can conclude
that for the systems we have constructed, functions $b_i$ in the
perturbation formula (HP) should satisfy two conditions:

\begin{itemize}
\item[(QH)]{\bf (quasi-homogeneity axiom)} {\it The obtained
system of Hamiltonian equations has to be quasi-homogeneous.}

In such a case, a Kowalevski matrix exists and we come to the last
condition. Suppose the invariant relations correspond to equations
$z_1=0,\dots, z_k=0$.

Denote by $p$ number of Casimirs: $n=p+2m$, where $2m$ is the
dimension of a general symplectic leaf.

\item[(ArA)]{\bf (Arithmetic axiom)} {\it For any nonzero solution
$C=(0,..., 0, c_{k+1},..., c_n)$ of the system \eqref{eq:Kowmat},
the Kowalevski matrix $K(C)$ has $n-p$  eigen-vectors tangent to
the symplectic leaf and $p$ transversal to it. Half of the
Kowalevski exponents which correspond to tangential eigen-vectors
and all of transversal ones are rational numbers. Irrational
numbers among the second half of tangential Kowalevski exponents
are divided into pairs such that the differences are integrally
dependent.}
\end{itemize}

The axioms of systems of Hess-Appel'rot type provide conditions
which determine classical Hess-Appel'rot system among
three-dimensional systems of Hess-Appel'rot type. More precisely,
suppose the two Poisson brackets are given on $e(3)$ as above and
a system is given by a Hamiltonian
\begin{equation}\label{eq:HALag}
H_{1}=H_0+JbM_3,
\end{equation}
where $H_0$ is the Hamiltonian of the Lagrange top corresponding to
the first Poisson structure, $M_3$ is its integral and a Casimir for
the second structure, $J$ is a nonzero constant and $b$ is a
function, such that the axioms of the systems of Hess-Appel'rot are
satisfied.

In \cite{DrGa2} the following {\it rigidity theorem} has been
proved.

\begin{thm}[\cite{DrGa2}] The  only non-zero polynomials $b$ which give systems of
Hess-Appel'rot type by relation \eqref{eq:HALag} are of the form
$$
b(z_1,\dots, z_6)=z_1 + kz_3.
$$
All systems of Hess-Appel'rot type of the form \eqref{eq:HALag}
are the classical Hess-Appel'rot systems.
\end{thm}

The last theorem provides a strong argument in favor of the choice
of axioms which have been postulated in \cite{DrGa2}. We take these
axioms as the starting point of our current research.

Since the Zhukovskii property is main object of our study in this
paper, we are going to focus ourselves on the first set of axioms.

\section{Examples of systems of Hess-Appel'rot
type}\label{sec:examples1}

\subsection{Rigid body systems on $so(n)\times
so(n)$}\label{subsec:examples1} The Euler-Poisson equations of
motion of a heavy rigid body fixed at a point are Hamiltonian on the
Lie algebra $e(3)$, which is the semi-direct product of Lie algebras
$\mathbb{R}^3$ and $so(3)$. Since $\mathbb{R}^3$ is isomorphic to
$so(3)$, there are two natural higher-dimensional generalizations of
Euler-Poisson equations. The first one is given by Ratiu in
\cite{Ra2} and it is to the semi-direct product $so(n)\times so(n)$
and the second one is to the Lie algebra $e(n)=\mathbb{R}^n\times
so(n)$.

Equations of a heavy $n$-dimensional rigid body on $so(n)\times
so(n)$ are :
\begin{equation}
\begin{aligned}
\dot M&=[M, \Omega]+[\Gamma, \chi]\\
\dot\Gamma&=[\Gamma,\Omega],
\end{aligned}
\label{3.1}
\end{equation}
where $M, \Omega, \Gamma,\chi\in so(n)$, $\Omega=AM$ and $\chi$ is
a constant matrix (see \cite{Ra2}). Here $A: so(n)\to so(n)$ is
the inverse of the rigid body kinetic energy operator.  The
Euler-Poisson equations \eqref{3.1} are Hamiltonian with the
Hamiltonian function
\begin{equation}
\label{hf-ha}
H=\frac12\langle
M,\Omega\rangle+\langle\chi,\Gamma\rangle=-\frac14\tr(M\Omega)-\frac12\tr(\chi\Gamma),
\end{equation}
in the standard Poisson structure on the semi-direct product
$so(n)\times so(n)$:
$$
\{M_{ij}, M_{jk}\}_1=-M_{ik},\qquad \{M_{ij},
\Gamma_{jk}\}_1=-\Gamma_{ik},\qquad
\{\Gamma_{ij},\Gamma_{kl}\}_1=0.
$$

The Casimir functions are $\tr(\Gamma^{2k})$, $\tr(M\Gamma^{2k+1})
$ and the dimension of generic symplectic leaf is
${n}(n-1)-2\left[\frac n2\right]$.

We will suppose that
\begin{equation}
\Omega=JM+MJ,\label{MANAKOV}
\end{equation}
where $J$ is a constant symmetric matrix. The operator $M \mapsto
JM+MJ$ belongs to the class of Manakov operators
 \cite{Ma} on $so(n)$.

The Lax representation together with Zhukovskii geometric
interpretation presented in the second section, where inspiration
for a construction of a higher-dimensional generalization of
Hess-Appel'rot system in \cite{DrGa2}. Let us first consider the
four-dimensional case.

\begin{dfn} The four-dimensional Hess-Appel'rot system
is described by the equations \eqref{3.1}, \eqref{MANAKOV}
and conditions:
\begin{equation}\label{ha-4-c}
J=\left (\begin{matrix} J_1&0&J_{13}&0\\
                    0&J_1&0&J_{24}\\
                    J_{13}&0&J_3&0\\
                    0&J_{24}&0&J_3
\end{matrix}\right), \qquad
\chi=\left(\begin{matrix} 0&\chi_{12}&0&0\\
                    -\chi_{12}&0&0&0\\
                    0&0&0&\chi_{34}\\
                    0&0&-\chi_{34}&0
                    \end{matrix}\right), \end{equation}
such that the operator $M\mapsto JM+MJ$ is positive definite and
$\chi_{12}^2+\chi_{34}^2\ne 0$.
\end{dfn}

The invariant surfaces are determined in the next lemma.

\begin{lem} {\rm (i)} For the four-dimensional Hess-Appel'rot system,
the following relations take place:
$$
\begin{aligned}
\dot M_{12}&=J_{13}(M_{13}M_{12}+M_{24}M_{34})+J_{24}(M_{13}M_{34}+M_{12}M_{24}),\\
\dot M_{34}&=J_{13}(-M_{13}M_{34}-M_{12}M_{24})+J_{24}(-M_{13}M_{12}-M_{24}M_{34}).
\end{aligned}
$$

{\rm (ii)} The system has two invariant relations:
\begin{equation}
M_{12}=0,\qquad M_{34}=0. \label{ha-4-ir}
\end{equation}
\end{lem}

Now let us introduce a new Poisson structure, compatible with the
standard one, as follows:
\begin{equation}
\begin{aligned}
&\{\Gamma_{ij},\Gamma_{jk}\}_2=-\Gamma_{ik},\quad
&\{M_{ij},\Gamma_{kl}\}_2=0,  \quad \{M_{13},
M_{23}\}_2=-\chi_{12},\\
&\{M_{14}, M_{24}\}_2=-\chi_{12}, \quad &\{M_{13},
M_{14}\}_2=-\chi_{34},\quad \{M_{23}, M_{24}\}_2=-\chi_{34}.
\end{aligned}
\label{ps1}
\end{equation}
Casimir functions in this structure are $M_{12}$, $M_{34}$,
$\Gamma_{12}^2+\Gamma_{13}^2+\Gamma_{14}^2+\Gamma_{23}^2+\Gamma_{24}^2+\Gamma_{34}^2$,
and
$\Gamma_{12}\Gamma_{34}+\Gamma_{23}\Gamma_{14}-\Gamma_{13}\Gamma_{24}$.

The situation with four-dimensional Hess-Appel'rot case is similar
to the three-dime\-nsional case: the  Hamiltonian is again a
quadratic deformation:
$$
H_{HA}=H_{LB}+J_{13}(-M_{12}M_{23}+M_{14}M_{34})+J_{24}(M_{12}M_{14}-M_{23}M_{34}),
$$
where
$H_{LB}=\frac12(2J_1M_{12}^2+(J_1+J_3)M_{13}^2+(J_1+J_3)M_{14}^2+
(J_1+J_3)M_{23}^2+(J_1+J_3)M_{24}^2+2J_3M_{34}^2)
+\chi_{12}\Gamma_{12}+\chi_{34}\Gamma_{34}$ is the Hamiltonian
function of the {\it Lagrange bitop}. The Lagrange bitop is a
complete integrable system of a heavy rigid body on $so(4)\times
so(4)$ defined in \cite{DrGa} and studied in details in
\cite{DrGa1}. Moreover, the Lagrange bitop is a bi-Hamiltonian
system. Assume that $J_1=a$, $J_3=1-a$. Then
$$
\begin{aligned} \tilde H_{LB}=&
\frac{(2a-1)(\chi_{12}M_{12}+\chi_{34}M_{34})}{\chi_{12}^2-\chi_{34}^2}
\left(\frac{M_{13}^2+ M_{14}^2+
M_{23}^2+M_{24}^2}2+\chi_{12}\Gamma_{12}+\chi_{34}\Gamma_{34}\right)\\
&+\frac{(1-2a)(\chi_{12}M_{34}+\chi_{34}M_{12})}{\chi_{12}^2-\chi_{34}^2}
(M_{23}M_{14}-M_{13}M_{24}+\chi_{12}\Gamma_{34}+\chi_{34}\Gamma_{12})\\
&+M_{12}\Gamma_{12}+M_{13}\Gamma_{13}+M_{14}\Gamma_{14}+M_{23}\Gamma_{23}+
M_{24}\Gamma_{24}+M_{34}\Gamma_{34}.
\end{aligned}
$$
is the Hamiltonian in the second structure \eqref{ps1}. The
functions $M_{12}$ and $M_{34}$, giving invariant relations for
the four-dimensional Hess-Appel'rot system, are Casimir
functions for the Poisson structure \eqref{ps1}.

Let us pass to the arbitrary dimension $n>4$.

\begin{dfn} The $n$-dimensional Hess-Appel'rot system
is described by the equations \eqref{3.1}, \eqref{MANAKOV},
together with conditions:
\begin{equation}\label{ha-n-c}
\begin{aligned}
&J=\diag(J_1,J_1,J_3,\dots,J_3)+J_{13}(E_1\otimes E_3+E_3\otimes
E_1)+J_{24}(E_2\otimes E_4+E_4\otimes E_2),
\\
&\chi=\chi_{12}E_1\wedge E_2, \qquad \chi_{12} \ne 0.
\end{aligned}
\end{equation}
\end{dfn}

Invariant relations are given in next lemma.

\begin{lem} {\rm (i)} For the $n$-dimensional Hess-Appel'rot system we have:
$$
\begin{aligned}
\dot M_{12}=&J_{13}(M_{12}M_{13}+M_{24}M_{34}+\sum_{p=5}^n
M_{2p}M_{3p})+\\
&J_{24}(M_{12}M_{24}+M_{13}M_{34}-\sum_{p=5}^nM_{1p}M_{4p})\\
\dot M_{34}=&-J_{13}(M_{13}M_{34}+M_{24}M_{12}+\sum_{p=5}^nM_{1p}
M_{p4})-\\& J_{24}(M_{13}M_{12}+M_{24}M_{34}+\sum_{p=5}^nM_{2p}
M_{3p})\\
\dot M_{3p}=&-J_{13}(M_{13}M_{3p}+M_{2p}M_{12})-J_{24}(M_{34}M_{2p}+M_{23}M_{4p})+\\
&M_{34}\Omega_{4p}-\Omega_{34}M_{4p}+\sum_{k=5}^n(M_{3k}\Omega_{kp}-\Omega_{3k}M_{4p}),\qquad
p>4,\\
\dot M_{4p}=&J_{13}(-M_{14}M_{3p}+M_{1p}M_{34})+J_{24}(M_{12}M_{1p}-M_{24}M_{4p})-\\
&M_{34}\Omega_{3p}+\Omega_{34}M_{3p}+\sum_{k=5}^n(M_{4k}\Omega_{kp}-\Omega_{4k}M_{4p}),\qquad
p>4,\\
\dot M_{kl}=&0,\qquad k,l>4
\end{aligned}
$$

{\rm (ii)} The system has the following set of invariant
relations:
\begin{equation}
M_{12}=0,\qquad M_{lp}=0,\qquad l,p\ge 3. \label{ha-n-ir}
\end{equation}
\end{lem}

In \cite{Ra2} the $n$-dimensional {\it Lagrange top}  on the
semidirect product $so(n)\times so(n)$ is constructed. In the metric $\Omega=JM+MJ$,
where $J=\diag(J_1,J_1,J_3,\dots,J_3)$, the $n$-dimensional Lagrange top is defined with Hamiltonian
$$
H_{L}=\frac12\left(2J_1M_{12}^2+(J_1+J_3)\sum_{p=3}^{n}(M_{1p}^2+M_{2p}^2)+2J_3\sum_{3\leq
p<q\leq n}M_{pq}^2\right)+\chi_{12}\Gamma_{12}
$$

Consider a Poisson structure
\begin{equation}
\{\Gamma_{ij},\Gamma_{jk}\}_2=-\Gamma_{ik},\,\{M_{ij},
M_{kl}\}_2=0,\,\,\{M_{1l}, M_{2l}\}_2=-1,\quad l=3,\dots,n,
\label{ps2}
\end{equation}
compatible with the standard one. The dimension of a symplectic
leaf in this structure is $\frac{(n-2)(n-3)}{2}
-\left[\frac{n-2}{2}\right]+4(n-2) $, hence there are
$\frac{n^2-5n+8}{2}+\left[\frac{n}{2}\right]$ Casimir functions:
$$
M_{12}, \quad M_{pq}, \quad \tr(\Gamma^{2k}), \qquad 2<p<q\le n,
\quad  k=1,\dots,\left[\frac{n}2\right].
$$

The $n$-dimensional Lagrange top is also bi-Hamiltonian system. In the
Poisson structure \eqref{ps2} its Hamiltonian is:
$$
\begin{aligned} \tilde
H_{L}&=(2a-1)M_{12}\left(\frac12\sum_{p=3}^n(M_{1p}^2+M_{2p}^2)+\Gamma_{12}\right)+\\
&(1-2a)\sum_{3\leq p<q\leq
n}M_{pq}(M_{1q}M_{2p}-M_{2q}M_{1p}+\Gamma_{pq})+\sum_{1\leq
p<q\leq n}M_{pq}\Gamma_{pq},
\end{aligned}
$$
where $J_{1}=a,\ J_{3}=1-a,\ \chi_{12}=1$.

Similarly as in dimension 3 and 4, Hamiltonian for the
Hess-Appel'rot system in arbitrary dimension $n$ is a quadratic
deformation of the Hamiltonian for the $n$-dimensional Lagrange
top:
$$
H_{HA}=H_{L}+\sum_{k=1}^n(J_{13}M_{1k}M_{3k}+J_{24}M_{2k}M_{4k}),
$$
and functions $M_{12}, M_{pq}, p,q\geq 3$, which give the
invariant relations,  are Casimir functions for the Poisson
structure \eqref{ps2}

Next theorem gives a Lax pair for the Hess-Appel'rot system.

\begin{thm} [\cite{DrGa2}] On invariant manifold given by the invariant relations,
the equations of $n$-dimensional Hess-Appel'rot system are equivalent to the
matrix equation
$$
\dot L(\lambda)=[L(\lambda), A(\lambda)],\qquad
L(\lambda)={\lambda^2}\frac1{J_1+J_3}\chi + \lambda
M+\Gamma,\qquad A(\lambda)=\lambda\chi +\Omega.
$$
\end{thm}

Using this Lax representation, the both classical and
algebro-geometric integration procedures are presented in
\cite{DrGa2} in dimension four.

\subsection{Rigid body systems on $e(n)$}
Let us now consider rigid body motion on the Lie algebra $e(n)$.
The standard Poisson structure on $e(n)^*\cong e(n)(M,\Gamma)$ is
given with:
$$
\{M_{ij},M_{kl}\}_1=-\delta_{jk}M_{il}, \qquad
\{\Gamma_i,\Gamma_j\}_1=0, \qquad \{M_{ij},\Gamma_k\}_1=-\Gamma_i
\delta_{jk}+\Gamma_j\delta_{ik}.
$$
Here we identified  $e(n)^*\cong e(n)$ by the use of a
non-invariant scalar product: $$ \langle
(M,\Gamma),(M,\Gamma)\rangle=\langle M,M\rangle+\langle
\Gamma,\Gamma\rangle=-\frac12\tr(MM)+\sum_{i=1}^n \Gamma_i^2. $$

A heavy rigid body Hamiltonian read: $H=\frac12\langle
M,\Omega\rangle+\langle\chi,\Gamma\rangle$, where $\chi\in\R^n$ is
the vector of mass centre and $\Gamma\in\R^n$ is a vertical vector
considered in the moving coordinate system and, as above,
$\Omega=AM$. We can choose $\chi=\chi_n E_n=(0,\dots,0,\chi_n)$.
The corresponding Euler--Poisson equations are:
\begin{equation}\label{EPC}
\begin{aligned}
\dot M &=[M,\Omega ]+\chi_n E_n \wedge \Gamma \\
\dot \Gamma &= -\Omega\cdot \Gamma.
\end{aligned}
\end{equation}

In \cite{Be} Belyaev considered the $n$-dimensional {\it Lagrange
top} defined by the Hamiltonian function:
\begin{equation}
H_\Lambda=\frac{a}{2}\langle M_\mathfrak d,M_\mathfrak d\rangle
+\frac{b}{2}\langle M_{\g},M_{\g}\rangle+\chi_n\Gamma_n,
\label{lagranz}
\end{equation}
where $so(n)=\g\oplus\D$ is orthogonal symmetric-pair
decomposition of $so(n)$: $ \g=\langle E_i \wedge E_j \, \vert\,
1\le i<j \le n-1 \rangle\cong so(n-1)$, $\D=\langle E_i \wedge E_n
\,\vert\, 1 \le i \le n-1\rangle$ and $a,b>0$ are real parameters.
Note that the kinetic energy has the Manakov form \eqref{MANAKOV},
where we take $\Omega=J_\Lambda M+MJ_\Lambda$,
$J_\Lambda=\diag(J_1,J_1,\dots,J_1,J_n)$.  Then $a=J_1+J_n$,
$b=2J_1$.

Belyaev proved noncommutative integrability of the system
\cite{Be}. The Lax representation is given by Reyman and
Semenov-Tian-Shanski \cite{RS}:
\begin{equation} \dot L(\lambda)=[L(\lambda),A(\lambda)],\quad
L(\lambda)=\hat\Gamma+\lambda  \hat M  + \lambda^2
\frac{\chi_n}{a} \hat E_n, \quad A(\lambda)=\hat\omega+ \lambda
\chi\hat E_n. \label{LA_pair}
\end{equation}
Here, for a given $(M,\Gamma)\in e(n)$, $\hat M, \hat\Gamma\in
so(n+1)$ are defined by
$$
\hat M= \left(
\begin{array}{cc}
M & 0  \\
0 & 0
\end{array} \right), \quad
\hat\Gamma= \left( \begin{array}{cc}
\mathbf{0} & \Gamma  \\
-\Gamma^t & 0
\end{array} \right).
$$

Let $\mathcal L$ be the set of integrals obtained from
\eqref{LA_pair}
\begin{equation}
\mathcal L: \qquad \tr(L(\lambda)^{2k}), \qquad k=1,\dots,\rank
SO(n+1), \qquad \lambda\in\R\label{LLL}
\end{equation}
and let $\mathcal S$ be linear functions on $\g=so(n-1)$
\begin{equation}
\mathcal S=\{M_{ij}\,\vert\, 1\le i<j\le n-1\}.
\label{S}\end{equation}
Then $\{\mathcal L,\mathcal L\}_1=0$,
$\{\mathcal L,\mathcal S\}_1=0$ and $\mathcal L+\mathcal S$ is
complete set of integrals of the Lagrange top system
 \eqref{EPC}, \eqref{lagranz}.

As above, let us define a Hess-Appel'rot system on $e(n)$ as a
perturbation of the Lagrange top using operator \eqref{MANAKOV} with
\begin{equation}
J=J_\Lambda+J_\Pi, \quad J_\Lambda=\diag(J_1,J_1,\dots,J_1,J_n),
\quad J_\Pi=J_{1n} (E_1\otimes E_n+E_n\otimes E_1). \label{ha-op}
\end{equation}

\begin{prop}
The equations \eqref{EPC}, \eqref{MANAKOV}, \eqref{ha-op} have
invariant relations:
\begin{equation}
M_{ij}=0, \qquad 1\le i<j \le n-1.\label{ha-inv-rel}
\end{equation}
Restriction of the system on invariant manifold is given with:
\begin{eqnarray}
&&\dot M_{in}=-\sum_{j=1}^{n-1} \Omega^\Pi_{ij} M_{jn} - \chi_n\,
\Gamma_i,
\quad i=1,\dots,n-1, \nonumber \\
&&\dot \Gamma_i=-(J_1+J_n)\Gamma_n M_{in}-\sum_{j=1}^{n-1}
\Omega^\Pi_{ij} \Gamma_j, \quad i=1,\dots,n-1, \label{Hess}\\
&&\dot \Gamma_n=(J_1+J_n) \sum_{j=1}^{n-1} \Gamma_j M_{jn},
\nonumber
\end{eqnarray}
where $ \Omega_\Pi=M J_\Pi+J_\Pi
M=-J_{1n}\sum_{i=2}^{n-1}(M_{in}E_1\wedge E_i+M_{1i}E_i\wedge
E_n). $
\end{prop}

\begin{rem} It can be proved that equations \eqref{EPC} have an invariant relation
 \eqref{ha-inv-rel} if and only if
$\pr_\D\circ A\circ\pr_\D$ is proportional to the identity
operator on $\D$ \cite{Jo2}.
\end{rem}

Applying a general construction of compatible Poisson structures
related to the symmetric pair decomposition of semi-simple Lie
algebras for the symmetric pair $(so(n+1),so(n))$ (see \cite{RS,
TF, Bo}), one gets the second Poisson structure on $e(n)$:
\begin{equation}\label{l-sps}
\{\Gamma_i,M_{jn}\}_2=\delta_{ij},\qquad
\{\Gamma_i,\Gamma_j\}_2=-M_{ij}.
\end{equation}

Let us denote $\chi_n=a=J_1+J_n=1$, $b=2J_1$.

\begin{prop}
{\rm (i)} The Casimirs of the structure \eqref{l-sps} are
functions \eqref{S} together with $
I=\frac12(M_{1n}^2+M_{2n}^2+\dots+M_{n-1,n}^2)+\Gamma_n. $

{\rm (ii)} The Lagrange top is Hamiltonian with respect to the
Poisson structure \eqref{l-sps} and the Hamiltonian function:
$$
H_{\Lambda,2}=\frac{1}{2}(\Gamma,\Gamma)+(b-1)\sum_{i,j=1}^{n-1}\Gamma_iM_{ij}M_{jn}\,.
$$
In this setting, the system is integrable in the commutative sense
by means of integrals \eqref{LLL}. The invariant tori are of
dimension $n-2$.
\end{prop}

Therefore, the system \eqref{EPC}, \eqref{MANAKOV}, \eqref{ha-op}
satisfies axioms
 (A1), (A2'), (HP) and (BP).

The restricted Hess-Appel'rot system (\ref{Hess}) admits the Lax
representation \eqref{LA_pair} with $M$ and $\Omega$ related by
\eqref{MANAKOV},  \eqref{ha-op}. But, only 3 integrals from the
family \eqref{LLL} are independent when invariant relations
\eqref{ha-inv-rel} are satisfied. Namely, the spectral curve is
then given by \cite{Jo2}:
\begin{eqnarray}
&& p(\lambda,\mu)=\det(L(\lambda)-\mu\mathrm{Id})=(-\mu)^{n-3}\left(\mu^4+\mu^2 P(\lambda)+Q(\lambda)\right)=0,  \label{CURVE}\\
&& P(\lambda)=F_2+\frac{2}{a}\lambda^2 F_1 + \lambda^4
\left(\frac{\chi_n}{a}\right)^2, \quad Q(\lambda)=\lambda^2
{F_3},\nonumber
\end{eqnarray}
where $F_1,F_2,F_3$ are integrals
\begin{eqnarray}
&&F_1=H=\frac{a}2 \sum_{i=1}^{n-1} M_{in}^2 +\chi_n\, \Gamma_n,
\quad F_2=\sum_{i=1}^{n} \Gamma_i^2=1, \nonumber \\
&&F_3={\sum_{1\le i<j \le n-1}(M_{in}\Gamma_j-M_{jn}\Gamma_i)^2},
\label{integrals2}
\end{eqnarray}
that correspond to the integrals \eqref{eq:intrigbody} of the
classical problem.

On the invariant submanifold (\ref{ha-inv-rel}), integrals $F_1,
F_2, F_3$ are not enough to give an integrability of the Lagrange
top. But, when conditions (\ref{ha-inv-rel}) are satisfied, the
system \eqref{EPC}, \eqref{lagranz} has additional integrals
$$
F_{ij}=M_{in}\Gamma_j-M_{jn}\Gamma_i, \qquad 1\le i<j \le n-1
$$
that implies that the Lagrange top on invariant submanifold
(\ref{ha-inv-rel}) can be solved by quadratures:
(\ref{ha-inv-rel}) is almost everywhere foliated on invariant
 2-dimensional tori.

Integrals $F_{ij}$ are not integrals of Hess-Appel'rot system
(\ref{Hess}). Nevertheless, the system \eqref{Hess} has additional
invariant relations
\begin{equation}
F_{ij}=M_{in}\Gamma_j-M_{jn}\Gamma_i=0, \qquad 1\le i<j \le n-1,
\label{dodatne}
\end{equation}
which are equivalent to condition $F_3=0$. Relations
\eqref{dodatne} show collinearity  of $\R^{n-1}$--vectors $
(M_{1n},M_{2n},\dots,M_{n-1,n})$ and $
(\Gamma_1,\Gamma_2,\dots,\Gamma_{n-1})$.

\subsection{Mishchenko--Fomenko flows}\label{mf-flows}
Let $G$ be a compact Lie group, $\mathfrak g$ its Lie algebra,
$\langle\cdot,\cdot\rangle$ a $\Ad_G$-invariant scalar product on
$\mathfrak g$. Let $a\in\mathfrak g$ be an arbitrary element
$\mathfrak g_a=\{ \eta\in \g, [a,\eta]=0\}$ be the isotropy
algebra and $G_a$ be the adjoint isotropy group of the element
$a$.

Let $\mathfrak g=\g_a\oplus\mathfrak d$ be the orthogonal
decomposition. Consider the linear operator (so called {\it
sectional operator} \cite{TF}) $A_{a,b,C}:\mathfrak g\to \mathfrak
g$, defined by
\begin{equation}
A_{a,b,C}(\xi)=\ad^{-1}_{a} \circ \ad_{b} \circ \pr_{\mathfrak d}
(\xi) + C (\pr_{\mathfrak g_a} \xi), \label{sectional}
\end{equation}
 where $b$ belongs  to the center of $\mathfrak g_a$,
 $\pr_{\mathfrak d}$
and $\pr_{\g_a}$  are the orthogonal (with respect to
$\langle\cdot,\cdot \rangle$) projections to $\mathfrak d$ and
$\mathfrak g_a$, respectively and $C: \mathfrak g_a \to \mathfrak
g_a$ is a positive definite, symmetric operator such that the
quadratic form $\langle \xi,C(\xi)\rangle$ is
$\Ad_{G_a}$-invariant:
\begin{equation}
[\xi, C (\xi)]=0, \quad \xi\in \g_a. \label{C}
\end{equation}

We can always find $b$ and $C$ such that $A_{a,b,C}$ is positive
definite. If $a$ is regular, i.e., $G_a$ is commutative, then the
condition (\ref{C}) is always satisfied.

Identify $\mathfrak g^*$ with $\mathfrak g$ by means of the scalar
product $\langle\cdot,\cdot\rangle$ and consider the
left-tri\-vi\-ali\-zation: $ T^*G \cong_l G\times \mathfrak
g=\{(g,\xi)\}. $ Then the quadratic form
\begin{equation}
H_{a,b,C}=\frac{1}{2}\langle A_{a,b,C}(\xi),\xi\rangle \label{hab}
\end{equation}
can be regarded as the Hamiltonian of a left-invariant Riemannian
metric $\kappa_{a,b,C}$ on $G$. After left $G$--reduction, the
equations of the geodesic flow take the form of Euler equations
\begin{equation}
\dot \xi=[\xi,\nabla H_{a,b,C}(\xi)]=[\xi,A_{a,b,C}(\xi)],
\label{b-flow}
\end{equation}
that are Hamiltonian with respect to the Lie-Poisson bracket
\begin{equation}
\{F_1,F_2\}_{1}=-\langle \xi, [\nabla F_1,\nabla
F_2]\rangle.\label{pb-1}
\end{equation}
The Casimir functions of the Lie-Poisson brackets
are invariant polynomials.

Let $\mathcal S$ be the set of linear functions on $\mathfrak
g_a$. Let $p_1,\dots,p_{\rank G}$ be the basic invariant
polynomials of the Lie algebra $\mathfrak g$. Mishchenko and
Fomenko proved that the polynomials obtained by shifting of
argument of invariants:
\begin{equation}
\mathcal C: \qquad p_{i,\lambda}(\xi)=p_i(\xi+\lambda a), \qquad
i=1,\dots,\rank  G, \qquad  \lambda\in\R \label{shifting}
\end{equation}
are commuting integrals of the system (\ref{b-flow}) and
$\{\mathcal C,\mathcal S\}_1=0$ \cite{MF1}. The algebra $\mathcal
C+\mathcal S$ is a complete algebra integrals of the system
(\ref{b-flow}) \cite{MF1, Bo, TF}. We refer to systems
(\ref{b-flow}) as {\it Mishchenko--Fomenko flows}.

Moreover, the system is bi-Hamiltonian (see \cite{Me, TF}). The
second Poisson structure is linear:
\begin{equation}
\{F_1,F_2\}_{2}=-\langle a, [\nabla F_1,\nabla F_2]\rangle,
\label{pb-2}\end{equation} with a set of Casimirs $\mathcal S$.
The integrals $p_{i,\lambda}(\xi)$, $i=1,\dots,\rank G$ are
Casimir functions of the bracket $
\{\cdot,\cdot\}_{1}+\lambda\{\cdot,\cdot\}_2$.

Now, let us perturb the Hamiltonian (\ref{hab}) as follows:
\begin{equation}
H_{a,b,C',D}=\frac{1}{2}\langle A_{a,b,C'}(\xi),\xi\rangle +
\langle D(\pr_{\mathfrak d} \xi),\xi \rangle,
\label{per-mf}\end{equation} where $D: \mathfrak d \to \g_a$ and
$A_{a,b,C'}$ is given by (\ref{sectional}) such that
$H_{a,b,C',D}$ is positive definite. Here $C'$ not need to satisfy
(\ref{C}).

\begin{prop}
The perturbed system
\begin{equation}
\dot \xi=[\xi,\nabla H_{a,b,C',D}(\xi)]
\label{c-flow}\end{equation} has an invariant manifold
\begin{equation}
\pr_{\mathfrak g_a}(\xi)=0 \label{c-zero}
\end{equation} and
satisfies axioms (A1), (A2), (HP), (BP) of the systems of the
Hess--Appelrot type.
\end{prop}

\section{Partial reductions}\label{sec:reduction}

\subsection{Classical Hess-Appel'rot system and spherical
pendulum}\label{ha-sp} Let us return to the classical problem
\eqref{eq:EPeq}, \eqref{Zhukovski}. Consider the motion of the
rigid body in the space reference frame. Let $R\in SO(3)$ be the
matrix that maps the moving reference frame to the fixed one.
Following Arnol'd's notation \cite{Ar}, denote
$\mathbf{m}=R\mathbf{M}, \mathbf{\omega}=R\mathbf{\Omega},
\mathbf{\gamma}=R\mathbf{\Gamma}$. The position of mass centre is
$z_0 \mathbf{e_3}$, where $\mathbf{e_3}=R\mathbf{E_3}=R(0,0,1)^T$.

Beside geometric interpretation of the conditions \eqref{eq:HAdef},
Zhukovski also noticed (see \cite{Zh}) that the motion of
$\mathbf{e_3}$ is described by the spherical pendulum equation.
Indeed, the Euler--Poisson equations for the variables
$(\mathbf{m},\mathbf{e_3})$, in the space frame, read
\begin{eqnarray}
&& \dot{\mathbf{m}}=z_0 \mathbf{\gamma}
\times\mathbf{e_3},\label{RHA1}\\
&& \dot{\mathbf{e_3}}=\omega \times
\mathbf{e_3}=\pr_{\mathbf{e_3^\perp}} \mathbf{\omega}\times
\mathbf{e_3}.\label{RHA2}
\end{eqnarray}

 On the invariant set $M_3=0$ we have
 $ \pr_{\mathbf{E_3^\perp}} \mathbf{\Omega}=(\Omega_1,\Omega_2,0)=(J_1 M_1,J_1 M_2,0)= J_1 \mathbf{M}$. The
 relation
 in the space frame gives
 \begin{equation}\label{RHA3}
\pr_{\mathbf{e_3^\perp}} \mathbf{\omega}= J_1 \mathbf{m}.
\end{equation}
By differentiation of \eqref{RHA2}, using \eqref{RHA1} and
\eqref{RHA3}, we get
\begin{equation*}
\ddot{\mathbf{e_3}}= J_1 z_0 (\gamma \times \mathbf{e_3})\times
\mathbf{e_3}+ J_1^2\mathbf{m}\times(\mathbf{m}\times
\mathbf{e_3})=-J_1z_0\mathbf{\gamma}+\mathbf{e_3}
(J_1z_0\langle\mathbf{e_3},\mathbf{\gamma}\rangle-\langle\dot{\mathbf{e_3}},\dot{\mathbf{e_3}}\rangle),
\end{equation*}
describing the Euler--Lagrange equations with multiplier
$\lambda=J_1z_0\langle\mathbf{e_3},\mathbf{\gamma}\rangle-\langle\dot{\mathbf{e_3}},\dot{\mathbf{e_3}}\rangle$
on the sphere $\langle{\mathbf{e_3}},{\mathbf{e_3}}\rangle=1$.
This is the pendulum system with Lagrangian
$$
l(\mathbf{e_3},\dot{\mathbf{e_3}})=\frac1{2J_1}\langle\dot{\mathbf{e_3}},\dot{\mathbf{e_3}}\rangle-z_0\langle\mathbf{\gamma},\mathbf{e_3}\rangle.
$$
In particular, the motion of $\mathbf{e_3}$ can be found by
elliptic quadratures.

\subsection{Invariant relations and reductions}
Note that, considered on the whole phase space $T^*SO(3)$ of the
rigid body motion, the function $M_3$ is the momentum mapping of
the right $SO(2)$-action-rotations of the body around the line
directed to the center of the mass. By the analogy with the
reduction of the system to the spherical pendulum, in this
subsection we study reductions of the Hamiltonian flows that
satisfy axioms (A1), (A2), (HP) of the systems of the
Hess--Appelrot type restricted to their invariant submanifolds.
Apparently, the lowering of order in Hamiltonian systems having
invariant relations was firstly studied by Levi-Civita (e.g., see
\cite[ch.~X]{LC}).

Let $G$ be a connected Lie group with a free proper Hamiltonian
action on a symplectic manifold $(M,\omega)$ with the momentum map
(\ref{moment_map}). Assume that $0$ is a regular value of $\Phi$.
Let $\xi_1,\dots,\xi_p$ be the base of $\g$. Then the zero level
set of the momentum mapping (\ref{moment_map}) is given by the
equations
\begin{equation} \label{moments}
M_0: \qquad \phi_i=(\Phi,\xi_i)=0, \qquad i=1,\dots,p.
\end{equation}

Let $(N_0,\omega_0)$ be the symplectic reduced space and let
$\pi_0: M_0 \to N_0=M_0/G$ be the canonical projection (see
subsection \ref{sim-red}). Consider the Hamiltonian equations:
\begin{equation} \label{hamiltonian}
\dot x=X_h(x).
\end{equation}

\begin{thm} [\cite{Jo2, Jo3}] \label{main}
\emph{(i)} Suppose that the restriction of $h$ to \eqref{moments}
is a $G$-invariant function. Then $M_0$ is an invariant manifold
of the Hamiltonian system \eqref{hamiltonian} and $X_{h}|_{M_0}$
projects to the Hamiltonian vector field $X_{h_0}$:
$d\pi_0(X_h)|_x=X_{h_0}|_{\pi_0(x)}$, where $h_0$ is the induced
function on $N_0$ defined by
\begin{equation} \label{induced_ham}
h|_{M_0}=\pi_0^* h_\eta=h_0\circ\pi_0.
\end{equation}

\emph{(ii)} The inverse statement also holds: if \eqref{moments}
is an invariant submanifold of the Hamiltonian system
\eqref{hamiltonian}, then the restriction of $h$ to $M_0$ is a
$G$-invariant function and $X_{h}|_{M_0}$ projects to the
Hamiltonian vector field $X_{h_0}$ on $N_0$, where $h_0$ is
defined by \eqref{induced_ham}.
\end{thm}

In both cases, the Hamiltonian vector field $X_h$ is not assumed
to be $G$-invariant on $M$. Moreover $X_h|_{M_0}$ may not be
$G$-invariant as well. It is invariant modulo the kernel of
$d\pi_0$, which is sufficient the tools of symplectic reduction
are still applicable.

We shall refer to the passing from $\dot x=X_h|_{M_0}$ to
\begin{equation} \label{REDUCED}
\dot y=X_{h_0}
\end{equation}
as a \emph{partial reduction}.

The natural way to obtain a Hamiltonian $h$ such that
corresponding flow have invariant relations (\ref{moments}) is to
perturb a $G$--invariant Hamiltonian $h_\Lambda\in C^\infty_G(M)$
\begin{equation}
h=h_\Lambda+\sum_{i=1}^p h_i\phi_i, \label{perturbation}
\end{equation}
where $h_i$ are arbitrary smooth functions on $M$. Then the set of
$p+1$ functions $h, \phi_1,\dots,\phi_k$ satisfy axioms (A1),
(A2'), (HP) of the systems of the Hess--Appelrot type, i.e., the
Hamiltonian system (\ref{hamiltonian}) is {\it restrectively
integrable} and of the form of {\it Hamiltonian perturbation} (see
\cite{DrGa2}).

The function $h_\Lambda$ is $G$--invariant and we can perform the
usual symplectic reduction to the Hamiltonian flow on the
symplectic reduced space $N_0$. Since the Hamiltonians $h$ and
$h_\Lambda$ coincide on $M_0$, the reduced flow is the same as the
partially reduced flow (\ref{REDUCED}). However, Hamiltonian
vector fields $X_h|_{M_0}$ and $X_{h_\Lambda}|_{M_0}$ are
different.

The partial reduction can be seen as a special case of the
symplectic reductions studied in \cite{BT, Li} (see also
\cite{LM}, Theorem 14.6, Ch. III).

\begin{thm} [\cite{LM}]
Let $(M,\omega,h)$ be a Hamiltonian system and let $M_0\subset M$
be an invariant submanifold of the vector field $X_h$ upon which
the symplectic form $\omega$ induces a 2-form $\omega_{M_0}$ of
constant rank. Assume that there is a surjective submersion with
connected fibres $\pi_0: M_0\to N_0$ onto another symplectic
manifold $(N_0,\omega_0)$, which satisfies
$\pi^*_0\omega_0=\omega|_{M_0}$. Then there exist a unique reduced
Hamiltonian function $h_0$ on $N_0$ such that
$h\vert_{M_0}=h_0\circ\pi_0$ and
$d\pi_0(X_h)\vert_x=X_{h_0}\vert_{\pi_0(x)}$.
\end{thm}

\subsection{Zhukovskii property}
An immediate corollary of theorem \ref{main} is

\begin{cor}\label{POSLED}
Suppose that the partially reduced system \eqref{REDUCED} is
completely integrable in the non-commutative sense and $N_0$ is
almost everywhere foliated on $r$-di\-me\-nsi\-onal isotropic
invariant manifolds, level sets of integrals $f^0_i$,
$i=1,\dots,\allowbreak\dim N_0-r$. Then $M_0$ is almost everywhere
foliated on $(r+\dim G)$-invariant isotropic manifolds
\begin{equation} \label{M}
\mathcal M_c=\{f_i=\pi^*_0 f^0_i=c_i \, \vert \, i=1,\dots, {\dim
N_0-r} \}
\end{equation}
of the system \eqref{hamiltonian}.
\end{cor}

In the case of the classical Hess--Appel'rot system we have
$M=T^*SO(3)$,  $G=SO(2)$ (rotations of the body around the vector
$\mathbf{E_3}$) and the reduced system is the spherical pendulum.
The reduced phase phase $N_0=T^*S^2$ is foliated on
two-dimensional tori, while the invariant manifold $M_0\subset
T^*SO(3)$ defined by $M_3=0$ is foliated on $3$-dimensional
invariant Lagrangian tori.

\begin{dfn}
We shall say that the Hamiltonian system (\ref{hamiltonian}) has
\emph{the Zhukovskii property} if it has invariant relations of the
form (\ref{moments}) and that the reduced system (\ref{REDUCED}) is
completely integrable.
\end{dfn}

Let us make  a terminological note:  in \cite{Jo2, Jo3} for {\it
the Zhukovskii property} has been used a different name, {\it
partial integrability} (or {\it geometrical Hess--Appol'rot
conditions} for natural mechanical systems). However, there is an
another notion named partial integrability which has been
introduced in a study based on the
Poincare-–Lyapo\-unov-–Nekhoroshev theorem \cite{N2, GG, N3} (see
also \cite{GMS, MPY}). There, a Hamiltonian system
(\ref{hamiltonian}) restricted to an invariant submanifold
$N\subset M$ of lower dimension is \emph{completely integrable},
i.e., $N$ is filled with periodic or quasi-periodic trajectories
of (\ref{hamiltonian}).

The following example will illustrate the difference between the
usual and partial reduction of Hamiltonian flows.

\begin{exm}\label{TORUS}
Let $G$ be a torus $\T^p$ and let the reduced flow be completely
integrable in the commutative sense by integrals
$f^0_1,\dots,f^0_m$, $m=\frac12\dim N_0$. Consider the regular
compact connected component level set of $f^0_1,\dots,f^0_m$. By
Liouville's theorem, it is diffeomorphic to a $m$-dimensional
torus $\T^m$ with quasi-periodic flow of (\ref{REDUCED}). Thus the
compact connected component $\hat{\mathcal M}_c=\pi^{-1}_0(\T^m)$
of (\ref{M}) is a \emph{torus bundle} over $\T^m$:
\begin{equation} \label{BUNDLE}
\begin{aligned}
\T^p\;\longrightarrow\;\;&\hat{\mathcal M}_c \\
&\;\Big\downarrow{}^{\pi_0} \\
&\;\T^m
\end{aligned}
\end{equation}

Suppose that $f_1,\dots,f_m$ can be extended to commuting $\T^p$-invariant functions
in some $\T^p$-invariant neighborhood $V$ of $\hat{\mathcal M}_c$.
Then, within $V$, $\hat{\mathcal M}_c$ is given by the equations
$$
f_1=c_1,\dots, f_m=c_m, \quad \phi_{1}=0,\dots,\phi_p=0.
$$

From the Noether theorem the functions $\phi_i$ commute with all
$\T^n$-invariant functions on $M$ and the following commuting
relations hold on $V$:
\begin{gather*}
\{f_a,f_b\}=\{f_a,\phi_i\}=\{\phi_i,\phi_j\}=0,\\
 a,b=1,\dots,m, \quad i,j=1,\dots,p.
\end{gather*}

Now, as in the case of commutative integrability of Hamiltonian
systems $\hat{\mathcal M}_c$ is a Lagrangian torus with tangent
space spanned by $X_{f_a}$, $X_{\phi_i}$, i.e., the bundle
(\ref{BUNDLE}) is trivial. However, in general, the flow of $X_h$
over the torus $\hat{\mathcal M}_c$ is {\it not quasi-periodic}:
the vector field $X_{h}$ does not commute with vector fields
$$
X_{f_1},\dots,X_{f_m},X_{\phi_1},\dots,X_{\phi_p}
$$
although Poisson brackets $\{h,f_a\}$, $\{h,\phi_i\}$ vanish on
$\hat{\mathcal M}_c$. So, in general, we can not apply the Lie
theorem \cite{Koz} to solve the system by quadratures.

Note that, if $h$ is $G$-invariant, the complete integrability of
the reduced and original system are closely related \cite{Zu,
Jo3}. In particular, if as above $G$ is torus, then the
reconstruction problem is easily solvable by quadratures (e.g.,
see \cite{MMR}). By contrary, in the case of partial reductions,
although the reduced motion $y(t)$ is found by quadratures, the
reconstruction problem in determining $x(t)$ ($\pi_0(x(t))=y(t)$,
$x(t_0)\in \pi_0^{-1}(y(t_0))$) in general, leads to
non-autonomous equations.
\end{exm}

\subsection{Reductions of additional symmetries}\label{additional-symmetries}
Suppose that an additional free Hamiltonian action of a connected
compact Lie group $K$ is given. Let
$$
\Psi: M \to \mathfrak k^*
$$
be the corresponding momentum mapping and let
$$
\sigma: M\to M/K
$$
be the natural projection. In what follows, by capital letters we
shall denote the functions on $P=M/K$ and with small letters
corresponding $K$-invariant functions on $M$:
$$
\sigma^*: C^\infty(P) \overset{\thickapprox}{\longrightarrow}
C^\infty_K(M), \qquad  F \mapsto f=F\circ \sigma.
$$

Since the action of $K$ is Hamiltonian, the Poisson bracket of two
$K$--invariant functions is $K$--invariant as well, so  manifold
$P=M/K$ carries the induced Poisson structure $\{\cdot,\cdot\}^K$
defined by $ \sigma^*\{F_1,F_2\}^K=\{\sigma^*F_1,\sigma^* F_2\}. $

The symplectic leaves in $(P, \{\cdot,\cdot\}^K)$ are of the form
$\Psi^{-1}(\mathcal O_\eta)/K$, where $\mathcal O_\eta$ is a
coadjoint orbit of $\eta\in\mathfrak k^*$. The Casimir functions
$I_1,\dots, I_r$  of $(P,\{\cdot,\cdot\}^K)$ are:
$$
I_j=(\sigma^*)^{-1}(p_j\circ\Psi), \qquad j=1,\dots,r=\rank K,
$$
where $p_1,\dots,p_r$ is the base of homogeneous invariants on
$\mathfrak k^*$, $r=\rank K$.

Further, suppose that the actions of $G$ and $K$ \emph{commute},
that is we have $ \{{\phi_i},{\psi_j}\}=0$, $i=1,\dots,p$,
$j=1,\dots,q$, where
\begin{equation*} \label{PSI}
\psi_j=(\Psi, \zeta_j), \qquad j=1,\dots,q.
\end{equation*}
and $\zeta_1,\dots,\zeta_q$ is a base of $\mathfrak k$. Thus
$\phi_i$ are $K$-invariant and we have induced $G$--action on $P$
given by the Hamiltonian functions $\Phi_i$, $
\Phi_i=(\sigma^*)^{-1}\phi_i$, $i=1,\dots,p.$

Suppose the functions $h_\Lambda,h_1,\dots,h_p$ in the definition
of the perturbed Hamiltonian (\ref{perturbation}) are
$K$-invariant functions. Then we can reduce the system
(\ref{hamiltonian}) to the Poisson manifold $P$ as well:
\begin{equation} \label{RF}
\dot F=\{F,H\}^K, \qquad F\in C^\infty(P).
\end{equation}
Here the reduced Hamiltonian $H$ is of the form
\begin{equation}
H=H_\Lambda+\sum_{i=1}^p H_i\Phi_i, \label{R-perturbation}
\end{equation}
where $H_\Lambda$ is $G$--invariant, i.e.,
$\{H_\Lambda,\Phi_i\}^K=0$. As a result we obtain reduced
invariant relations
\begin{equation} \label{RIM}
\Phi_i=0, \quad i=1,\dots,p
\end{equation}
defining the invariant manifold $P_0=M_0/K$ of the reduced flow
(\ref{RF}).

The functions $\Phi_i$ are integrals of the non-perturbed flow so
$H$ satisfy axioms (A1), (A2') and (HP) of systems of Hess-Appel'rot
type.

On the other hand, we have the induced Hamiltonian $K$-action on
$(N_0,\omega_0)$ with the momentum mapping $\Psi_0$ satisfying
$\Psi|_{M_0}=\Psi_0 \circ \pi_0$ (see diagram \eqref{diagram}
below). Since $h$ is $K$-invariant, the reduced Hamiltonian $h_0$ on
$N_0$ is also $K$-invariant.  Therefore the momentum mapping
$\Psi_0$ is conserved along the flow of the partially reduced system
(\ref{REDUCED}).
\begin{align}
\mathfrak k^* \,\,\,&\,\,\,\, \qquad\mathfrak k^*
\qquad\qquad\qquad\quad\,\,\,\,\,
\mathfrak k^*\nonumber\\
_\Psi\Big\uparrow \;\;\;&\qquad _\Psi\Big\uparrow \qquad\qquad\qquad\quad _{\Psi_0}\Big\uparrow\nonumber\\
\mathfrak g^* \overset{\Phi}{\longleftarrow} M\;\;
&\supset\;\;M_0=\Phi^{-1}(0)\; \overset{\pi_0}{\longrightarrow} N_0=M_0/G\label{diagram}\\
^\sigma\Big\downarrow \;\;\;&\qquad ^\sigma\Big\downarrow \nonumber\\
P=M/K &\supset \,\,P_0=M_0/K\nonumber
\end{align}

\begin{rem}\label{primedba}
Suppose we have additional integrals of the reduced system
(\ref{REDUCED}) implying Zhukovskii property of
(\ref{hamiltonian}) with respect to the $G$--action. Then $M_0$ is
almost everywhere foliated on invariant isotropic manifolds
$\mathcal M_c$ (see Corollary \ref{POSLED}). Hence we obtain
invariant foliation of $P_0=M_0/K$ on manifolds of the form
$\mathcal P_c=(K\cdot \mathcal M_c)/K$. It is clear that $\mathcal
P_c$ are isotropic submanifolds of the appropriate symplectic
leaves in $(P,\{\cdot,\cdot\}^K)$.
\end{rem}

\subsection{Natural mechanical systems}\label{nms}
Let $(Q,\kappa,v)$ be a  natural mechanical system, where the metric
$\kappa$ is also regarded as a mapping $\kappa: TQ\to T^*Q$. Let $G$
be a connected Lie group acting freely and properly on $Q$ and
$\rho: Q\to Q/G$ be the canonical projection.

For Lagrangian systems, it is convenient to work with tangent
bundle reductions. Let $\mathcal V_q=\{\xi_q \; \vert \;
\xi\in\g\}$ be the tangent space to the fiber $G \cdot q$ ({\it
vertical space at} $q$) and $\mathcal V=\cup_q \mathcal V_q$ be
the vertical distribution. Consider the {\it horizontal
distribution} $\mathcal H$ orthogonal to $\mathcal V$ with respect
to the metric $\kappa$:
\begin{equation}
\mathcal H=\{X_q \in T_q Q\, \vert \, (\kappa_q (X_q) ,\,
\xi_q)=0, \,\, \xi\in\g,\,\, q\in Q\} \label{horizontal*}
\end{equation}

Let $\Phi$ be the cotangent bundle momentum mapping
\eqref{ctg-mm}. Since $\mathcal H=\kappa^{-1}(\Phi^{-1}(0))$, the
horizontal distribution is invariant with respect to the
``twisted" $G$-action
\begin{equation} \label{action}
{g}\diamond(q,X) =(g\cdot q,\,\kappa_{g\cdot
q}^{-1}\circ(dg^{-1})^*\circ\kappa_q(X)),\quad X\in T_q Q,
\end{equation}
that is the pull-back of canonical symplectic $G$-action on $T^*Q$
via metric $\kappa$. With the above notation we get \cite{Jo2}

\begin{thm}
\emph{(i) \ (Partial Noether theorem)} The horizontal distribution
\eqref{horizontal*} is an invariant submanifold of the
Euler--Lagrange equations \eqref{Lagrange} of the natural
mechanical system $(Q,\kappa,v)$ if and only if the potential $v$
and the restriction $\kappa_\mathcal H$ of the metric $\kappa$ to
$\mathcal H$ are $G$-invariant with respect to the action
\eqref{action}.

\emph{(ii) \ (Partial Lagrange--Routh reduction)} If $\mathcal H$
is an invariant submanifold of the system $(Q,\kappa,v)$, then the
trajectories $q(t)$ with velocities $\dot q(t)$ that belong to
$\mathcal H$ project to the trajectories $\pi(q(t)$ of the natural
mechanical system $(Q/G,K,V)$ with the potential $V(\pi(q))=v(q)$
and the metric $K$ obtained from $\kappa_\mathcal H$ via
identification $\mathcal H/G\cong T(Q/G)$.
\end{thm}

Note that when $\kappa$ is $G$-invariant, the twisted $G$-action
(\ref{action}) coincides with usual $G$-action:
$g\cdot(q,X)=(g\cdot q,dg(X))$ and the induced metric $K$ is the
submersion metric.

\section{Systems of  Hess--Appel'rot type and Zhukovskii
property}\label{sec:partial}

In the sequel we shall study the relationship between the
Zhukovskii property and systems of  Hess--Appel'rot type on the
Poisson manifold $(P,\{\cdot,\cdot\}^K)$. In particular, we shall
establish the dimension of invariant manifolds stated in Remark
\ref{primedba}.

We say that some property holds for a {\it generic point} $x$ of
the manifold $M$ if the property holds on an open dense set
$U\subset M$.

\subsection{Hamiltonian perturbation}
Consider the system (\ref{RF}) with Hamiltonian function $H$ of
the form (\ref{R-perturbation}), where $H_i$ are arbitrary smooth
functions on $P$ and $H_\Lambda\in C^\infty_G(P)$ is a
$G$--invariant Hamiltonian.

\begin{lem}\label{ocuvanje}
If $F$ is a $G$--invariant integral of the system
\begin{equation}
\dot F=\{F,H_\Lambda\}^K, \qquad F\in C^\infty(P),
\label{non-perturbed}
\end{equation}
then it is also a integral of the system \eqref{RF} with perturbed
Hamiltonian function \eqref{R-perturbation}, restricted to the
invariant manifold \eqref{RIM}.
\end{lem}

\begin{proof}
Directly, from $\{\Phi_i,F\}^K=0$ and $\{H_\Lambda, F\}^K=0$ we
get
\begin{eqnarray*}
\{H,F\}^K\vert_{P_0} &=& \{H_\Lambda,F\}^K\vert_{P_0}+\sum_{i=1}^p \{H_i\Phi_i,F\}^K\vert_{P_0}\\
&=& \sum_{i=1}^p
H_i\{\Phi_i,F_j\}^K\vert_{P_0}+\sum_{i=1}^p\Phi_i\{H_i,F\}^K\vert_{P_0}=0.
\end{eqnarray*}
\end{proof}

Therefore, if the non-perturbed flow (\ref{non-perturbed}) is
integrable, for a generic perturbation (\ref{perturbation}), only
$G$--invariant integrals remain to be integrals of the perturbed
flow restricted to $P_0$.

In what follows, we suppose that at a generic point $z$ of $P_0$,
functions $I_1,\dots,I_r$ are {independent}, or equivalently, the
symplectic leaf through $z$ is {regular}.

Let $F_1,\dots,F_\rho$ be the Poisson--commuting integrals of the
non-perturbed system (\ref{non-perturbed}), mutually independent
and independent of functions $I_1,\dots,I_r$,
$\Phi_1,\dots,\Phi_p$ at a generic point $z\in P_0$. Let
$$
s=\dim G_z=\dim (\langle dI_1,\dots,dI_r\rangle \cap \langle
d\Phi_1,\dots,d\Phi_p\rangle)\vert_z\,,
$$
at a generic $z\in P_0$ ($G_z\subset G$ is the isotropy group of
$z$ with respect to the induced $G$--action) and let
\begin{eqnarray}
&&\mathcal
F=\{\Phi_1,\dots,\Phi_p,I_1,\dots,I_r,F_1,\dots,F_\rho\}\subset C^\infty (P),\label{f}\\
&&\mathcal
F_0=\{\psi^0_1,\dots,\psi^0_q,f_1^0,\dots,f_\rho^0\}\subset
C^\infty(N_0).\label{f-0}
\end{eqnarray}

Here $\psi^0_i=(\Psi_0,\zeta_i)$ are component of the momentum
mapping $\Psi_0$ and $f_i^0$ are obtained from $F_i$ by the
composition
\begin{equation}
C^\infty_G(P) \overset{\thickapprox}{\longrightarrow}
C^\infty_{G,K}(M) \overset{\imath^*}{\longrightarrow}
C^\infty_K(M_0)\overset{\thickapprox}{\longrightarrow}
C^\infty_K(N_0),\label{composition}
\end{equation}
where $\imath: M_0 \hookrightarrow M$ is the inclusion.

\begin{thm}\label{procena}
\emph{(i)} \begin{equation} 2\rho \le \dim P+2s-2\dim G-\rank
K.\label{VII}
\end{equation}

 \emph{(ii)} If
$\mathcal F$ is a complete set at a generic point $z\in P_0$ then
$G$ is commutative and the partially reduced flow \eqref{REDUCED}
is completely integrable by means of integrals $\mathcal F_0$.

\emph{(iii)} Contrary, suppose $\mathcal F_0$ is a complete set of
functions on $N_0$. Then the manifold $P_0$ is almost everywhere
foliated on $\frac12(\dim P-\rank K)$--dimensional level sets of
integrals $\mathcal F$, both of the non-perturbed
\eqref{non-perturbed} and the perturbed system \eqref{RF}. This is
enough for integrability of the non-perturbed system
\eqref{non-perturbed} in the case when $G$ is a commutative Lie
group.
\end{thm}

\begin{proof} {\sc I}
Consider Theorem \ref{Bol_Jov} where we put $(M,\omega,G,\Phi)$ to
be equal to $(N_0,\omega_0,K,\Psi_0)$. We have (see \cite{BJ2}):
\begin{eqnarray}
&&\dind \Psi_0^*C^\infty(\mathfrak k^*)=\dind
C^\infty_K(N_0)=\nonumber\\
&& \qquad \qquad \qquad\,\,\,\,\,=\dind(\Psi^*C^\infty(\mathfrak
k^*)+C^\infty_K(N_0))=\dim K_\mu-\dim K_y\label{I}\\
&&\ddim C^\infty_K(N_0)=\dim N_0+\dim K_y-\dim K \nonumber
\end{eqnarray}
for a generic $y\in N_0$, $\mu=\Psi_0(y)$.

As above, let $p_1,\dots,p_r$ be the base of homogeneous
invariants on $\mathfrak k^*$, $r=\rank K$. Then functions
$$
i^0_i=p_i\circ \Psi_0, \qquad i=1,\dots,r
$$
belong to $\Psi^*C^\infty(\mathfrak k^*) \cap C^\infty_K(N_0)$,
and therefore belong to the centers in both algebras
$\Psi^*C^\infty(\mathfrak k^*)$ and  $C^\infty_K(N_0)$. Among them
there are exactly $\dind C^\infty_K(N_0)$ independent ones. Under
the composition (\ref{composition}), $i_j^0$ corresponds to the
Casimir function $I_j$ of the bracket $\{\cdot,\cdot\}^K$.

Since $F_1,\dots,F_\rho$ are mutually independent and independent
of the functions $I_1,\dots,I_r$ and $\Phi_1,\dots,\Phi_p$ at a
generic point $z\in P_0$, it follows that $f_1^0,\dots,f^0_\rho$
are also mutually independent and independent of the functions
$i_1^0,\dots,i_r^0$ and
\begin{equation}
\ddim \{i_1^0,\cdots,i^0_r,f_1^0,\dots,f_\rho^0\}=\dind
C^\infty_K(N_0)+\rho.\label{III}
\end{equation}

The set $\{i_1^0,\cdots,i_r^0,f_1^0,\dots,f_\rho^0\}$ is a
commutative subset of $C^\infty_K(N_0)$, so
\begin{equation}
\ddim\{i_1^0,\cdots,i_r^0,f_1^0,\dots,f_\rho^0\}\le \frac12(\ddim
C^\infty_K(N_0)+\dind C^\infty_K(N_0)) \label{IV}
\end{equation}

The condition that a symplectic leaf through $z\in P_0$ is regular
implies that the coadjoint orbit through $\mu=\Psi_0(y)$ is
regular, where $y=\pi_0(x)$, $z=\sigma(x)$, $x\in M_0$. Therefore,
from (\ref{I}) we get
\begin{eqnarray}
\ddim C^\infty_K(N_0)+\dind C^\infty_K(N_0)
&=& \dim N_0+\rank K-\dim K \nonumber\\
&=& \dim P-2\dim G+\rank K \label{V}
\end{eqnarray}

Also $s=\dim G_z=\dim K_y$, for a generic $y=\pi_0(x)$,
$z=\sigma(x)$, $x\in M_0$. Thus
\begin{equation}
\dind C^\infty_K(N_0)=\rank K-s\,.\label{VI}
\end{equation}

Finally, combining (\ref{III}), (\ref{IV}), (\ref{V}) and
(\ref{VI}) we obtain inequality (\ref{VII}).

{\sc II} By using the proof of item (i), we get
\begin{eqnarray}
&&\ddim \mathcal F=\dim G+\rank K+\rho-s, \nonumber\\
&&\dind\mathcal F=\rank G+\rank K+\rho-s, \nonumber \\
&&\ddim\mathcal F+\dind\mathcal F=\dim G+\rank G +2\rank
K+2\rho-2s\nonumber \\
&&\qquad   \le \dim G+\rank G +2\rank K-2s+\dim P+2s-2\dim G-\rank K\label{VIII}\\
&&\qquad   = \dim P+\rank G-\dim G+\rank K  \le \dim P+\rank
K.\label{IX}
\end{eqnarray}

If \eqref{f} is a complete set at a generic point $z\in P_0$ then
we have equalities both in (\ref{VIII}) and (\ref{IX}). Equality
in (\ref{IX}) implies that $G$ is a commutative group. Equality in
(\ref{VIII}) means that we have equality in (\ref{VII}), that is,
according (\ref{IV}),
$\{i_1^0,\cdots,i_r^0,f_1^0,\dots,f_\rho^0\}$ is a complete
commutative subset in $C^\infty_K(N_0)$.

Then, from Theorem \ref{Bol_Jov} we get that $
\{\psi^0_1,\dots,\psi^0_q,i_1^0,\cdots,i_r^0,f_1^0,\dots,f_\rho^0\}
$ is a complete set on $N_0$. Since $i_j^0$ are polynomial
functions in $\psi^0_i$, the set \eqref{f-0} will be also a
complete set of integrals of the partially reduced system
(\ref{REDUCED}).

{\sc III} Suppose \eqref{f-0} is a complete set. Then
$\{i_1^0,\cdots,i_r^0,f_1^0,\dots,f_\rho^0\}$ is a complete
commutative subset in $C^\infty_K(N_0)$ and we have equality in
\eqref{VII}. Therefore $ \ddim\mathcal F=\dim G+\rank
K+\frac12(\dim P+2s-2\dim G-\rank K)-s= \frac12(\dim P+\rank K)$.
The dimension of invariant level-sets given by $\mathcal F$ is
$\dim P-\ddim\mathcal F=\frac12(\dim P-\rank K)$.
\end{proof}

\subsection{Lifting of bi--Poisson structure}
As above, consider the system
(\ref{RF}) with Ha\-miltonian function $H$ of the form
(\ref{R-perturbation}), where $H_i$ are arbitrary smooth functions
on $P$ and $H_\Lambda\in C^\infty_G(P)$. Suppose in addition to
the Poisson structure $\{\cdot,\cdot\}_1=\{\cdot,\cdot\}^K$, that
the non-perturbed system (\ref{non-perturbed}) is Hamiltonian with
respect to the another Poisson structure $\{\cdot,\cdot\}_2$ which
is compatible with the first one. Also, we suppose that functions
$\Phi_i$ are Casimir functions of the second bracket:
\begin{equation}
\{\Phi_i,F\}_2\equiv 0, \qquad i=1,\dots,p, \quad F\in
C^\infty(P). \label{casimir}
\end{equation}

Thus, the Hamiltonian flow (\ref{RF}) satisfies axioms (A1), (A2'),
(HP) and (BP) of systems of  Hess-Appel'rot type.

Let
$$
\Pi=\{\{\cdot,\cdot\}_{\lambda_1,\lambda_2}=\lambda_1\{\cdot,\cdot\}_1+\lambda_2\{\cdot,\cdot\}_2\;
\vert\; \lambda_1,\lambda_2\in \mathbb{R},
\;\lambda_1^2+\lambda_2^2\ne 0\},
$$
be the corresponding pencil of compatible Poisson structures.

\begin{thm}\label{prenosenje}
The pencil of compatible Poisson structures $\Pi$ on $P$ induces
the pencil of compatible Poisson structures $\Pi_0$ within the
algebra $C^\infty_K(N_0)$ of $K$--invariant functions on $N_0$
\end{thm}

\begin{proof}
The Poisson brackets $\{\cdot,\cdot\}_1$ and $\{\cdot,\cdot\}_2$
are compatible. That is why it is enough to prove that both
Poisson brackets are correctly defined within the algebra of
functions $C^\infty_K(N_0)$.

The $G$ and $K$ actions on $(M,\omega)$ is Hamiltonian. Thus
$\{\cdot,\cdot\}_1$ is well defined in all of the algebras:
$C^\infty_K(M)\cong C^\infty (P)$, $C^\infty(M)$, $C^\infty_G(M)$,
$C^\infty(N_0)$ and $C^\infty_K(N_0)$.

The function $F\in C^\infty(P)$ is $G$--invariant if and only if
\begin{equation}
\{F,\Phi_i\}_1=0, \qquad i=1,\dots,p. \label{G-inv}
\end{equation}

The bracket $\{\cdot,\cdot\}_2$ induces the Poisson structure
within $C^\infty_G(P)$ if the bracket $\{F_1,F_2\}_2$ of two
$G$--invariant functions is again a $G$--invariant function.

Suppose $F_1$ and $F_2$ are $G$--invariant. Let us write the
Jacobi identity for the bracket
$\{\cdot,\cdot\}_1+\{\cdot,\cdot\}_2$ and functions
$F_1,F_2,\Phi_i$:
\begin{eqnarray*}
&&\{\{F_1,F_2\}_1+\{F_1,F_2\}_2,\Phi_i\}_1+\{\{F_1,F_2\}_1+\{F_1,F_2\}_2,\Phi_i\}_2+\\
&&\{\{\Phi_i,F_1\}_1+\{\Phi_i,F_1\}_2,F_2\}_1+\{\{\Phi_i,F_1\}_1+\{\Phi_i,F_1\}_2,F_2\}_2+\\
&&\{\{F_2,\Phi_i\}_1+\{F_2,\Phi_i\}_2,F_1\}_1+\{\{F_2,\Phi_i\}_1+\{F_2,\Phi_i\}_2,F_1\}_2=0.
\end{eqnarray*}

Taking into account identities (\ref{casimir}) and (\ref{G-inv})
for $\Phi_i$, $F_1, F_2$, from the Jacobi identity we get
\begin{equation}\label{G-inv*}
\{\{F_1,F_2\}_1,\Phi_i\}_1+\{\{F_1,F_2\}_2,\Phi_i\}_1=0.
\end{equation}
Since $F_1$ and $F_2$ are $G$--invariant, $\{F_1,F_2\}_1$ is also
$G$--invariant, so the first term is equal to zero. Hence, the
second term in (\ref{G-inv*}) equals zero as well and the bracket
$\{\cdot,\cdot\}_2$ induces the Poisson structure within
$C^\infty_G(P)\cong C^\infty_{K,G}(M)$.

The functions $\Phi_i$ are Casimir functions for the bracket
$\{\cdot,\cdot\}_2$. Whence, after the restriction, we have well
defined bracket $\{\cdot,\cdot\}_2$ on the submanifold
(\ref{RIM}). From the above considerations it follows that the
Poisson bracket $\{\cdot,\cdot\}_2$ is well defined within the
algebra of functions $C^\infty_G(P_0)\cong
C^\infty_{G,K}(M_0)\cong C^\infty_K(N_0)$.
\end{proof}

The pencil $\Pi_0$ can be used as a tool in proving complete
integrability of the partially reduced system (\ref{REDUCED}). For
example,  suppose that all structures in $\Pi$, except possible
$\{\cdot,\cdot\}_2$, have the maximal rank equal to $\dim P-\rank
K$, and that their Casimir functions are globally defined on $P$.
Then the union $\mathcal C$ of Casimir functions of all Poisson
structures from $\Pi$ non-proportional to $\{\cdot,\cdot\}_2$ is a
commutative set with respect to the all Poisson structures from
$\Pi$. Moreover, let $\mathcal S$ be the set of Casimir functions
of the bracket $\{\cdot,\cdot\}_2$. Then $\{\mathcal C,\mathcal
S\}_1=0$ and, since $\{\Phi_1,\dots,\Phi_p\}\subset \mathcal S$,
the functions in $\mathcal C$ are $G$--invariant. If the pencil
$\Pi$ satisfies conditions of the form \eqref{uslovi-bolsinova}
(for more details see \cite{Bo}), the set of functions $\mathcal
F=\mathcal C+\mathcal S$ will be a complete set of functions.
However, as we have seen in Theorem \ref{procena}, even if
$\mathcal F$ is complete, it not need to be complete at the points
of the invariant set $P_0$.

Nevertheless, suppose that the Casimir functions $\mathcal C$
induce the Casimir functions $\mathcal C_0\subset C^\infty_K(N_0)$
of the brackets non-proportional to $\{\cdot,\cdot\}_2$. Then,
from Theorems \ref{bolsinov*} and \ref{Bol_Jov} we get

\begin{prop}\label{uslovi}
The set $\mathcal C_0$ is a complete commutative subset of
$C^\infty_K(N_0)$  if and only if the corank of the coplexified
bivectors
$$
\pi^\mathbb{C}_y: (\mathrm{ann}(K\cdot y))^\mathbb{C} \times
(\mathrm{ann}(K\cdot y))^\mathbb{C}\to \mathbb{C}
$$
is equal to $\dind C^\infty_K(N_0)$ for all
$\pi\in\Pi_0^\mathbb{C}$ at the generic point $y\in N_0$. In this
case the partially reduced system \eqref{REDUCED} is integrable by
means of integrals $\mathcal F_0=\mathcal
C_0+\{\psi_1^0,\dots,\psi^0_q\}.$
\end{prop}

Here $\mathrm{ann}(K\cdot y)\subset T^*_y N_0$ is the annihilator
of the tangent space to the orbit $K\cdot y$ of $y$.

The condition stated in Proposition \ref{uslovi} is verified in
proving the complete commutative integrability of geodesic flows
of normal metrics on adjoint orbits in \cite{BJ3, MP} (see Remark
\ref{compatible} given below).

\section{Zhukovskii property. Examples}\label{sec:examples2}

\subsection{Magnetic flows on adjoint orbits}
In order to describe partial reduction of the Hess--Appelrot rigid
body system \eqref{3.1}, \eqref{MANAKOV}, \eqref{ha-4-c},
\eqref{ha-n-c} and the geodesic flows (\ref{c-flow}) we shall need
a description of certain natural mechanical systems on adjoint
orbits recently studied, in the presence of the additional
magnetic force, in \cite{BJ4, BJ5}. Note that integrable magnetic
flows on homogeneous spaces are also given in \cite{RS, Ef, MSY}.

Let $G$ be a compact connected Lie group with the Lie algebra
$\mathfrak g$ and invariant scalar product $\langle \cdot,\cdot
\rangle$. Consider the $G$--adjoint orbit $\mathcal
O(a)\subset\mathfrak g$. The tangent space at $x=\Ad_g(a)$ is
simply the orthogonal complement to $ \mathfrak
g_x=\{\xi\in\mathfrak g\,\vert\,[x,\xi]=0\}. $ The cotangent
bundle $T^*\mathcal O(a)$ can be represented as a submanifold of
$\mathfrak g\times \mathfrak g$: $$ T^*\mathcal O(a)=\{(x,p)\,
\vert\, x=\Ad_g(a), p\in\mathfrak g_x^\perp\},
$$
with the paring between $p\in T^*_x \mathcal O(a) \cong \mathfrak
g_x^\perp$ and $\eta\in T_x\mathcal O(a)$ given by
$p(\eta)=\langle p,\eta\rangle$. Then the canonical symplectic
form $\omega$ on $T^*\mathcal O(a)$ can be seen as a restriction
of the canonical linear symplectic form of the ambient space
$\mathfrak g\times\mathfrak g$: $\sum_{i=1}^{\dim\mathfrak g} dp_i
\wedge dx_i$, where $p_i$, $x_i$ are coordinates of $p$ and $x$
with respect to some base of $\mathfrak g$. Let  $\Omega$ be the
standard Kirillov-Kostant symplectic form on $\mathcal O(a)$.

The canonical $G$--action $ g\cdot (x,p)=(\Ad_g x,\Ad_g p) $ on
the magnetic cotangent bundle $(T^*\mathcal
O(a),\omega+\epsilon\rho^*\Omega)$ is Hamiltonian with the
momentum mapping (see \cite{Ef, BJ4})
\begin{equation}
m: T^*\mathcal O(a) \to \g^*\cong \g, \quad m(x,p)=[x,p]+\epsilon
x. \label{mmm}
\end{equation}

A natural generalization of the magnetic spherical pendulum to the
orbit $\mathcal O(a)$ is the mechanical system with the kinetic
energy given by the {\it normal metric}  and the potential
function $V(x)=\langle c,x\rangle$, i.e., with the Hamiltonian
$$
h_c(x,p)=\frac12 \langle [x,p],[x,p]\rangle+ \langle c,x\rangle.
$$

Another natural class of systems are the magnetic geodesic flows
of the $G$--invariant metrics $K_{a,b}$ defined by the Hamiltonian
function
$$
h_{a,b}(x,p)=\frac12\langle [{b_x} p],[x,
p]\rangle=-\frac12\langle \ad_x\ad_{b_x}p,p\rangle,
\label{red_ham}
$$
where $b$ belongs to the center of $\mathfrak g_a$ and
$b_x=\Ad_g(b)$ for $x=\Ad_g(a)$ \cite{BJ4}. For compact groups, we
can take $b$ such that $K_{a,b}$ is positive definite. If $b=a$ we
get the Hamiltonian of the {normal metric}.

The equations of the magnetic pendulum, in redundant variables
$(x,p)$, are given by (see \cite{BJ5})
\begin{equation}\label{xx}
\begin{aligned}
& \dot x=[x,[p,x]],  \\
& \dot p=[p,[p,x]] +  \epsilon [x,p]-c+\pr_{\g_x} c,\end{aligned}
\end{equation}
while the magnetic geodesic flow read (see \cite{BJ4})
\begin{equation}\label{adjoint1}
\begin{aligned}
& \dot x=-\ad_x\ad_{b_x} p=[[b_x,p],x],  \\
& \dot p=-\ad_x^{-1}[p,[x,[b_x,p]]]+\pr_{\g_x}[[b_x,p],p]+\epsilon
[b_x,p]. \end{aligned}
\end{equation}

Let $\mathcal O(a)$ be an arbitrary orbit.

\begin{thm}[\cite{BJ4, BJ5}] \label{complete2}
The magnetic pendulum system, for a regular element $c\in\mathfrak
g$, and the geodesic flows of the metrics $K_{a,b}$ on
$(T^*\mathcal O(a),\omega+\epsilon\rho^*\Omega)$, described by
equations \eqref{xx} and \eqref{adjoint1},  respectively, are
completely integrable by means of polynomial integrals.
\end{thm}

\subsection{Partial reduction of  rigid body
systems} Consider the construction given is Section 4 for the case
when the symplectic manifold $M$ is the phase space of the
$n$-dimensional rigid body motion about a fixed point $T^*SO(n)$.
As usual, we use the left trivialization $T^*SO(n)\cong_l SO(n)
\times so(n)(g,M)$, where $so(n)$ and $so(n)^*$ are identified by
the use of invariant scalar product $\langle
X,Y\rangle=-\frac12\tr(XY)$.

We follow Ratiu's generalization of the heavy rigid body motion
\cite{Ra2}. Let
$$
\gamma=\gamma_{12} E_1 \wedge E_2+\gamma_{34} E_3\wedge E_4 +
\dots + \gamma_{k-1,k} E_{k-1}\wedge E_{k}, \quad k=2[n/2]=2\rank
SO(n),
$$
where $\gamma_{i,i+1}$ are mutually different. Also, let
\begin{align*}
&\chi=\chi_{12} E_1 \wedge E_2,  &n>4,\,\, \\
&\chi=\chi_{12} E_1 \wedge E_2+ \chi_{34} E_3 \wedge E_4,
&\,\,\,\,n=4,
\end{align*}
where $\chi_{12}\ne \chi_{34}$, $\chi_{12}\ne 0$. Then the adjoint
isotropy groups of $\gamma$ and $\chi$ and corresponding isotropy
Lie algebras read, respectively (e.g., see \cite{Besse})
\begin{align*}
&K=SO(n)_\gamma=\overset{k\,times}{\overbrace{SO(2)\times \cdots
\times SO(2)}}, &\mathfrak k=\langle E_1 \wedge
E_2,\dots,E_{k-1}\wedge E_k\rangle, \\
&G=SO(n)_\chi=SO(2)\times SO(n-2),  &\mathfrak g=\langle E_1
\wedge E_2, E_i \wedge E_j\,\vert\, 3 \le i<j \le n\rangle.
\end{align*}

Note that $K$ is a maximal torus in $SO(n)$, so the adjoint orbit
through $\gamma$ is the flag manifold $\mathcal
O(\gamma)=SO(n)/K$, while the adjoint orbit $\mathcal
O(\chi)=SO(n)/G=SO(n)/SO(2)\times SO(n-2)$ is Grassmannian variety
$Gr^+(n,2)$ of oriented 2-dimensi\-onal planes through the origin
in $\R^n$.

Consider the natural left action of $K$ and right action of $G$ on
$T^*SO(n)$. The corresponding momentum mapping $\Psi$ and $\Phi$,
in the left trivialization are given by
$$
\Psi(g,M)=\pr_{\mathfrak k}(\Ad_g M), \qquad
\Phi(g,M)=\pr_\mathfrak g(M)
$$

The Hamiltonian of the $n$--dimensional Hess--Appelrot rigid body
system is the Hamiltonian perturbation of the Lagrange system (or
Lagrange bi-top system for $n=4$):
\begin{equation}
h(g,M)=\frac12\langle M,\Omega\rangle+\langle \Ad_{g^{-1}}
\gamma,\chi\rangle, \label{HAH}
\end{equation} where $\Omega=JM+MJ$
and $J$ is given by \eqref{ha-n-c}. The fixed element $\gamma$
play the role of the horizontal vector (direction of the
gravitational force) as seen in the space reference frame.

\begin{lem}
The Hamiltonian \eqref{HAH} is left $K$--invariant and right
$G$--invariant on the zero level set of the momentum mapping
$\Phi$:
\begin{equation}
(T^*SO(n))_0: \,\, \phi_{12}(M,g)=M_{12}=0, \,\,
\phi_{ij}(M,g)=M_{ij}=0, \,\, 3 \le i <j \le n.\label{IHA}
\end{equation}
\end{lem}

We have $P=(T^*SO(n))/K \cong so(n)\times \mathcal
O(\gamma)(M,\Gamma)$. The $K$--reduced  flow is described by
equations \eqref{3.1}, \eqref{MANAKOV}, \eqref{ha-n-c} (or
\eqref{ha-4-c} for $n=4$), where one have to fix values of
invariants in $\Gamma$ in order that $\Gamma$ belongs to the
adjoint orbit $\mathcal O(\gamma)$

On the second hand, the symplectic reduced space
$N_0=(T^*SO(n))_0/G$ is symplectomorphic to the cotangent bundle
of the adjoint orbit $\mathcal O(\chi)\cong Gr^+(n,2)$.

Let $x=\Ad_g \chi\in \mathcal O(\chi)$ (in $3$--dimensional case,
$x$ represents the position of the mass center in the space
coordinates). We can rewrite the Hamiltonian (\ref{HAH}) in the
form
$$
h(M,g)=\frac{J_1+J_n}{2}\langle M,M\rangle + \phi_{12} H_{12} +
\sum_{3 \le i<j\le n} \phi_{ij} H_{ij} + \langle \gamma,x\rangle.
$$
Therefore
$$
h(M,g)\vert_{(T^*SO(n))_0}=\frac{J_1+J_n}{2}\langle M,M\rangle +
\langle \gamma,x\rangle. $$

After reduction to $T^*\mathcal O(\chi)$, the bi-invariant kinetic
energy term  goes to the kinetic energy of the normal metric
multiplied by $(J_1+J_n)$ (e.g, see
 \cite{Besse}). Thus, likewise in the 3--dimensional case,
 the partially reduced flow is the
pendulum  system on $\mathcal O(\chi)$:
\begin{equation}\label{p1}
\begin{aligned}
& \dot x=(J_1+J_n)[[x,p],x], \\
& \dot p=(J_1+J_n) [[x,p],p] -\gamma + \pr_{so(n)_x}
\gamma,\end{aligned}
\end{equation}
with Hamiltonian $ h_0(x,p)=\frac{J_1+J_n}{2}\langle [x,p],[x,p]
\rangle +\langle x,\gamma \rangle$ (we follow the notation of the
previous section). It follows from Theorem \ref{complete2} that
the reduced system is completely integrable. Hence the system
satisfies Zhukovskii property and we get the following qualitative
behavior of the system,

\begin{thm} \label{DG}
The partial reduction of the Hess-Appel'rot rigid body problem
defined by the Hamiltonian \eqref{HAH} is completely integrable
pendulum type system \eqref{p1} on the oriented Grassmannian
variety $Gr^+(n,2)$.  The invariant manifold \eqref{IHA} is almost
everywhere foliated by invariant $\dim SO(n)$-dimensional
Lagrangian invariant manifolds that project to the
$2(n-2)$-dimensional Liouville tori of the reduced system
\eqref{p1}.
\end{thm}

So, in this approach, after solving the pendulum type system, the
equations of Hess-Appel'rot rigid body problem reduces to $(\dim
SO(n-2)+1)$--differential equations of the reconstruction problem.

Similarly as in subsection \ref{ha-sp},  it can be proved that the
partial reduction of the rigid body system \eqref{EPC},
\eqref{MANAKOV}, (\ref{ha-op}) is the pendulum system on the
$(n-1)$--dimensional sphere $S^{n-1}$ (see \cite{Jo2}), given by
the Hamiltonian function
$$
h_0(x,p)=\frac{J_1+J_n}{2}(p,p)+\chi_n\langle x,\gamma\rangle.
$$

Here the cotangent bundle of the sphere is realized  as a
submanifold $(x,x)=1$, $(x,p)=0$ of $\R^{2n}(x,p)$ and $x,\gamma
\in S^{n-1}$ represent the direction of the position of the mass
center $\chi$ and the position of the vertical axes $\Gamma$ in
the space coordinates, respectively.

\subsection{Mishchenko--Fomenko flows}
We follow the notation of subsections \ref{mf-flows} and
\ref{additional-symmetries}, where we take
$(M,G,K,P,N_0)=(T^*G,G_a,G,\mathfrak g,T^*(G/G_a))$. Here we
consider the right $G_a$--action and the left $G$--action on
$T^*G$. The momentum maps read $\Phi(g,\xi)=\pr_{\mathfrak
g_a}(\xi)$ and $\Psi(g,\xi)=\Ad_g\xi$, respectively.

The metric $\kappa_{a,b,c}$ is right $G_a$--invariant and we can
project it to the homogeneous space $G/G_a$, that is to the
adjoint orbit of $a$. Since we deal with the right action, the
vertical distribution is left-invariant: $\mathcal V_g=g\cdot
\mathfrak g_a$, while from the definition of $\kappa_{a,b,c}$, the
horizontal distribution is $\mathcal H_g=g\cdot\mathfrak d$ and
the submersion metric does not depend on $c$. The submersion
metric is exactly the metric $K_{a,b}$ on $\mathcal O(a)$ defined
above.

The Hamiltonian $H_{a,b,C',D}$ defines left--invariant metric on
$G$ that we shall denote by $\kappa_{a,b,C',D}$. The Hamiltonian
functions $H_{a,b,C',D}$ and $H_{a,b,C}$ coincides on the
invariant manifold
\begin{equation}
(T^*G)_0 \cong_l G\times \mathfrak d\,. \label{b-zero}
\end{equation}

Therefore, the partial reduction of the geodesic flow of the
metric $\kappa_{a,b,C',D}$ is the geodesic flow of the metric
$K_{a,b}$. The flow (\ref{adjoint1}) is completely integrable for
any $\epsilon$. Thus the geodesic flow of the perturbed
Mishchenko--Fomenko metric $\kappa_{a,b,C,D}$ satisfies Zhukovskii
property.

\begin{rem}
If $a$ is a {\it regular} element of the Lie algebra $\g$ ($G_a$
is Abelian) then  the invariant manifold (\ref{b-zero}) is almost
everywhere foliated on invariant isotropic tori of the geodesic
flow $\kappa_{a,b,C',D}$ and the motion over the tori is not
quasi-periodic (see Example \ref{TORUS}). The same holds for the
left $G$--reduced flow (\ref{c-flow}) restricted to
(\ref{c-zero}). Namely, it can be proved that $\mathcal S \subset
\mathcal C$ and that $\mathcal C$ is complete at a generic point
$\xi$ that belongs to (\ref{c-zero}) (e.g, see Theorem 3.2 in
\cite{BJ3}). According Lemma \ref{ocuvanje}, the perturbed flow
(\ref{c-flow}) has the same foliation of (\ref{c-zero}) on
invariant tori as the non-perturbed flow (\ref{b-flow}).  Item
(ii) of Theorem \ref{procena} then implies the integrability of
geodesic flow of the metric $K_{a,b}$.
\end{rem}

\begin{rem}\label{compatible}
For a {\it singular} $a\in\g$, $\mathcal C+\mathcal S$ restricted
to (\ref{c-zero}) is {\it not complete}. The complete
integrability of the geodesic flow of the normal metric and the
magnetic geodesic flows (\ref{adjoint1}) on the adjoint orbit
$\mathcal O(a)$ follows from the completeness of the commutative
set $\mathcal C_0$ induced from \eqref{shifting} within the
algebra $C^\infty_G(T^*\mathcal O(a))$. The completeness is
obtained by verifying the condition stated in Proposition
\ref{uslovi} for the pencil of compatible Poisson structures
within $C^\infty_G(T^*\mathcal O(a))$, induced from (\ref{pb-1})
and (\ref{pb-2}), see \cite{BJ3, MP, BJ4}. Besides, item (iii) of
Theorem \ref{procena} gives us an estimate of the dimension of
invariant manifolds within (\ref{c-zero}) of the perturbed flow
(\ref{c-flow}).
\end{rem}

\begin{rem}
The horizontal distributions $\mathcal H'$ and $\mathcal H$ (see
subsection \ref{nms}) of the metrics $\kappa_{a,b,C,D}$ and
$\kappa_{a,b,C}$ are different
$$
\mathcal H'_g=\kappa_{a,b,C',D}^{-1}(g\cdot \mathfrak d)\ne g\cdot
\mathfrak d=\mathcal H_g\,,
$$
for $D\ne 0$, while $\mathcal H'_g=\mathcal H_g$ for $D=0$. Let
$D=0$. In this case, the integrals $\mathcal C$ remains to be the
integrals of the perturbed system (\ref{c-flow}) not only on the
invariant manifold (\ref{c-zero})) but on $\mathfrak g$ as well.
\end{rem}

\subsection{Singular Manakov flows}
Similar perturbations as those of the Mishchenko--Fomenko flows
can be performed for other integrable Euler equations with
symmetries.  The natural candidate is the singular Manakov flow.
Let $$
a=(\overset{k_1\,times}{\overbrace{\alpha_1,\dots,\alpha_1}},
\dots,\overset{k_r\,times}{\overbrace{\alpha_r,\dots,\alpha_r}}),
\quad b=(\overset{k_1\,times}{\overbrace{\beta_1,\dots,\beta_1}},
\dots,\overset{k_r\,times}{\overbrace{\beta_r,\dots,\beta_r}}), $$
where $k_1+k_2+\dots+k_r=n$, $\alpha_i \ne \alpha_j$, $\beta_i\ne
\beta_j$, $i,j=1,\dots,r$ and let
\begin{equation}so(n)=\g \oplus\mathfrak d=so(k_1)\oplus
so(k_2)\oplus\dots\oplus so(k_r)\oplus \mathfrak d \label{v}
\end{equation}
be the orthogonal decomposition, where $ \g=\{X\in so(n)\, \vert\,
[X,a]=0\}. $ By $M_{\g}$ and $M_\mathfrak d$ we denote the
projections of $M\in so(n)$ with respect to (\ref{v}). Further,
let $C: \g\to \g$ be an arbitrary positive definite operator. We
take $a$ and $b$ such that the sectional operator $A_{a,b,C}:
so(n)\to so(n)$ defined via
\begin{equation}
A_{a,b,C}(M_{\mathfrak d}+M_{\g})=\ad_a^{-1}\ad_b(M_{\mathfrak
d})+C(M_{\g}), \label{A}
\end{equation}
is positive definite. Here $\ad_a$ and $\ad_b$ are considered as
invertible linear transformations from $\mathfrak d$ to
$[a,\mathfrak d]\subset Sym(n)$. Let $H_{a,b,C}=\frac12\langle
M,A_{a,b,C} (M)\rangle$.
 We refer to Euler equations
\begin{equation}
\dot M=[M,\Omega], \quad \Omega= \nabla H_{a,b,C}(M)=A_{a,b,C}(M),
\label{Euler}
\end{equation}
as the {\it singular Manakov Flow}. The operator $A_{a,b,C}$
satisfies Manakov condition $[M,b]=[\Omega,a]$, so we have the Lax
representation with rational parameter $\lambda$ (see Manakov
\cite{Ma}):
\begin{equation}
\dot L(\lambda)=[L(\lambda),U(\lambda)], \quad
L(\lambda)=M+\lambda a, \quad U(\lambda)=\Omega+\lambda b.
\label{Lax}
\end{equation}

In the case the eigenvalues of $a$ are all distinct, i.e., $\g=0$,
Manakov proved  that the solutions of the Euler equations
(\ref{Euler}) are expressible in terms of $\theta$-functions by
using the algebro-geometric integration procedure developed by
Dubrovin in \cite{Du} (see \cite{Ma}).   The explicit verification
that integrals arising from the Lax representation
\begin{equation}
\mathcal L=\{\tr(M+\lambda a)^k\, \vert\, k=1,2,\dots,n,\,
\lambda\in \mathbb{R}\}, \label{integrals}
\end{equation}
 form a
complete commutative set on $so(n)$ was given by Mishchenko and
Fomenko \cite{MF1}.

Let us denote the set of linear functions on $\g$ by $\mathcal S$.
These are additional integrals in the case the eigenvalues of $C$
are not all distinct and $C$ is proportional to the identity
operator, or more generally in the case $C$ is an
$\Ad_{G}$-invariant, where $ G=SO(k_1)\times SO(k_2)\times
\dots\times SO(k_r)\subset SO(n). $ The complete integrability of
the system is proved by Bolsinov by using the pencil of compatible
Poisson brackets given by the canonical Lie-Poisson bivector
$$
\pi_1(\xi_1,\xi_2)\vert_M=-\langle M, [\xi_1,\xi_2]\rangle
$$
and
\begin{equation*}
\pi_2(\xi_1,\xi_2)\vert_M=-\langle M, \xi_1 A \xi_2-\xi_2 A \xi_1
\rangle\label{**}
\end{equation*}
(see \cite{Bo} and \cite{TF}, pages 241-244). The another proof is
given in \cite{DGJ}. Namely, we have $\{\mathcal L,\mathcal
S\}_{so(n)}=0$ and $\mathcal L+\mathcal S$ is complete at a
generic $M\in so(n)$.

Now, the perturbation follows the perturbation of
Mishchenko--Fomenko flows:
\begin{equation}
\dot M=[M,\Omega], \quad \Omega=\nabla H_{a,b,C,D}(M),
\label{manakov-perturbed}
\end{equation}
where $D:\mathfrak d \to \mathfrak g$ and
$H_{a,b,C,D}=H_{a,b,C}+\langle M_\mathfrak g,D(M_\mathfrak
d)\rangle$ is positive definite (we do not suppose that $C$ is
$\Ad_G$--invariant). The system \eqref{manakov-perturbed} has the
invariant manifold
\begin{equation}
\mathfrak d: \qquad M_{\mathfrak g}=0.\label{inv-man}
\end{equation}
Besides, the restriction of the system to \eqref{inv-man} has the
Manakov L-A pair \eqref{Lax} and integrals \eqref{integrals}.

Consider diagram \eqref{diagram}, where we take $$
(M,G,K,P,N_0)=(T^*SO(n),G,SO(n),so(n),T^*(SO(n)/G)), $$ with the
right $G$--action and the left $SO(n)$--action on $T^*SO(n)$. By
the use of the map \eqref{composition}, the commutative set of
function $\mathcal L$ induces a complete commutative set $\mathcal
L_0$ within $C^\infty_{SO(n)}(T^*(SO(n)/G))$ \cite{DGJ}. On the
other side, if $G$ is not commutative then $\mathcal L+\mathcal S$
is not complete at $\mathfrak d$.

\begin{prop} If $G$ is commutative, i.e., $\alpha_i \le 2$,
$i=1,\dots,r$, then the set of function $\mathcal L+\mathcal S$ is
a complete set at a generic point $M\in\mathfrak d$.
\end{prop}

\begin{proof} The proposition directly follows from item (iii) of Theorem \ref{procena}.
Alternatively, let $L_M=\{\nabla_M\tr(M+\lambda A)^k\, \vert\,
k=1,2,\dots,n,\, \lambda\in \mathbb{R}\}. $ According to
(\ref{coisotropic}), $\mathcal L+\mathcal S$ is complete at $M$ if
\begin{equation}
(L_M+\g)^{\pi_1}\subset L_M+\g. \label{LL}\end{equation}

The relation \eqref{LL} is proved in \cite{DGJ} by using Theorem
\ref{bolsinov}, namely by verifying that the dimension of the
linear spaces (25) and (26) in \cite{DGJ} are equal to $n$. The
dimension of the space (26) calculated in Lemma 2 \cite{DGJ} holds
for a generic $M\in\mathfrak d$. On the other side, the dimension
of (25) is equal to $n$ for elements in a generic position with
the property that $M_\mathfrak g$ is a regular element in
$\mathfrak d$ (see \cite{TF}, pages 234-237). Since $\g$ is
commutative, the dimension of (25) will be $n$ as required.
 \end{proof}

\begin{exm}
As an example, take $n=2r$,
$a=(\alpha_1,\alpha_1,\dots,\alpha_r,\alpha_r)$. Then
$\g=so(2)\otimes so(2)\dots \otimes so(2)$ is the Cartan
subalgebra. The set of integrals $\mathcal L+\mathcal S$ is
complete at $\mathfrak d$ so the invariant set $\mathfrak d$ of
the systems \eqref{Euler} and \eqref{manakov-perturbed} is
foliated on invariant tori. The matrix $L(\lambda)$ satisfies
$$
L_{12}=L_{21}=L_{34}=L_{43}=\dots=L_{2r-1,2r}=L_{2r,2r-1}=0.
$$
In other words, the systems \eqref{Euler}, restricted to
\eqref{inv-man}, is an example of integrable isoholomorhic system
(\cite{DrGa1}, see Remark \ref{isohol} given below).
\end{exm}

The bi-Hamiltonian formulation of singular Manakov flows can be
performed by using the pencil
$\Pi=\{\pi_{\lambda_1,\lambda_2}=\lambda_1\pi_1+\lambda_2\pi_2\}$
given above (see \cite{Bo, TF}). The singular brackets within the
pencil $\Pi$ are proportional to
$\pi_{1,-\alpha_1},\dots,\pi_{1,-\alpha_r}$ and the linear
function on $so(k_i)\subset \g$ are among the Casimirs of the
brackets $\pi_{1,-\alpha_i}$, $i=1,\dots,r$. So, if only one $k_i$
is greater of 1 (say $k_1=k>1$, $k_2=\dots=k_r=1$, $r=n-k+1$,
i.e., $\g=so(k)$), the perturbed singular Manakov flow
\eqref{manakov-perturbed} satisfies axiom (BP) with respect to the
second Poisson structure $\pi_{1,-\alpha_1}$. The partially
reduced system is completely integrable geodesic flow on the
Stiefel variety $SO(n)/SO(k)$.

\section{Integration of the magnetic pendulum on
$Gr^+(4,2)$}\label{sec:pendulum}

Let us consider closely the pendulum system given by the equations
\eqref{p1} in dimension four with the magnetic term added. Since
the orbit $\mathcal{O}(\chi)\cong Gr^+(4,2)$ is defined with
invariants, the cotangent bundle of the Grassmannian $Gr^+(4,2)$
is given by the constraints:
\begin{equation}
\label{t*o}
\begin{aligned}
&x_{12}^2+x_{13}^2+x_{14}^2+x_{23}^2+x_{24}^2+x_{34}^2=\chi_{12}^2+\chi_{34}^2,\\
&x_{34}x_{12}+x_{23}x_{14}-x_{13}x_{24}=\chi_{12}\chi_{34},\\
&x_{12}p_{12}
+x_{13}p_{13}+x_{14}p_{14}+x_{23}p_{23}+x_{24}p_{24}+x_{34}p_{34}=0,\\
&x_{34}p_{12}+x_{23}p_{14}-x_{13}p_{24}+p_{34}x_{12}+p_{23}x_{14}-p_{13}x_{24}=0.
\end{aligned}
\end{equation}

Introducing the magnetic momentum mapping \eqref{mmm}, the equations
\begin{equation*}
\begin{aligned}
& \dot x=(J_1+J_n)[x,[p,x]],  \\
& \dot p=(J_1+J_n)[p,[p,x]] +  \epsilon
[x,p]-\gamma+\pr_{so(4)_x}\gamma\end{aligned}
\end{equation*}
become:
\begin{equation}
\begin{aligned}
\dot m& =[\gamma, x]\\
\dot x&=(J_1+J_n)[m,x].
\end{aligned}
\label{8.1}
\end{equation}

The equations \eqref{8.1} are special case of the equations of the
completely symmetric Lagrange bitop. The Lagrange bitop is
defined in \cite{DrGa} and studied in details in \cite{DrGa1}.
Thus, the integration procedures given in \cite{DrGa, DrGa1}  can
be applied to the considered system. We will present here  both of
them, the classical and the algebro-geometric integration
procedures.

After solving the system \eqref{8.1} one has $m$ and $x$ as  known
functions of time. In order to find $p$ as a function of time, one
needs to solve the equation
$$m=[x, p]+\epsilon x$$
in $p$. Let us
 denote
 $$m_0 = m-\epsilon x=[x,p].$$
  Since $\langle [p,x],so(4)_x\rangle=\langle [x, so(4)_x], p\rangle=0$
 we have $m_0\in so(4)_x^{\perp}$. The operator $\ad_x:so(4)_x^\perp\mapsto so(4)_x^\perp$
 is bijective, thus a solution $p=\ad^{-1}_x(m_0)$ is unique.

\subsection{Classical integration of the magnetic pendulum}

Starting from a well- known decomposition  $so(4)=so(3)\oplus
so(3)$, let us introduce as in \cite{DrGa1}
$$
m_1=\frac 12 (m_++m_-),\qquad m_2=\frac 12 (m_+-m_-)
$$
 where $m_+, m_-$ are two
three-dimensional vectors which correspond to  four-dimensional matrix $m_{ij}$ according to
$$
(m_+,m_-)\longmapsto \left(\begin{matrix}
0 & -m^3_{+} & m^2_{+} & -m^1_{-}\\
m^3_{+} & 0 & -m^1_{+} & -m^2_{-}\\
-m^2_{+} & m^1_{+} & 0 & -m^3_{-}\\
m^1_{-} & m^2_{-} & m^3_{-} & 0
\end{matrix}\right).
$$
(Similar decomposition can be performed for $x, \gamma$).

Equations \eqref{8.1} become
\begin{equation}
\begin{aligned}
\dot m_i&=2(\gamma_i\times x_i),\qquad \dot
x_i=2(J_1+J_n)(m_i\times x_i),\qquad i=1,2,
\end{aligned}
\label{8.2}
\end{equation}
where
$$
\gamma_1=(0,0,-\frac12 (\gamma_{12}+\gamma_{34})),\quad
\gamma_2=(0,0,-\frac12 (\gamma_{12}-\gamma_{34})).
$$

If we denote $m_1=(p_1, q_1, r_1),\ m_2=(p_2, q_2, r_2)$,
then the first group of the equations \eqref{8.2} becomes
$$
\begin{aligned}
&{\dot p}_1=-2\gamma_{(1)3}x_{(1)2}, &&{\dot p}_2=-2\gamma_{(2)3}x_{(2)2},\\
&{\dot q}_1=2\gamma_{(1)3}x_{(1)1},  &&{\dot q}_2=2\gamma_{(2)3}x_{(2)1},\\
&{\dot r}_1=0, &&{\dot r}_2=0,
\end{aligned}
$$
where we denoted with $x_{(i)j}$ the $j$ component of the vector $x_i$.

The integrals of motion are for $i=1,2$:
$$
\begin{aligned}
&r_i=f_{i1}\\
&(J_1+J_n)(p_i^2+q_i^2)+2\gamma_{(i)3}x_{(i)3}=f_{i2}\\
&p_ix_{(i)1}+q_ix_{(i)2}+r_ix_{(i)3}=f_{i3}=\epsilon \chi_{(i)3}\\
&x_{(i)1}^2+x_{(i)2}^2+x_{(i)3}^2=f_{i4}=\chi_{(i)3}^2,
\end{aligned}
$$
The constants $f_{i3}$ and $f_{i4}$ are found from the conditions
\eqref{t*o}.

Following \cite{DrGa1} and introducing $\rho_i, \sigma_i$, defined
with $p_i=\rho_i\cos\sigma_i$, $q_i=\rho_i\sin\sigma_i$, we get
$$
\begin{aligned}
&{\dot\rho_i}^2+\rho_i^2{\dot\sigma}_i^2=4\gamma_{(i)3}^2(\chi_{(i)3}^2-x_{(i)3}^2),\\
&\rho_i^2\dot{\sigma}_i=2\gamma_{(i)3}(\epsilon\chi_{(i)3}-f_{i1}x_{(i)3}).
\end{aligned}
$$
It follows that $u_i=\rho_i^2$ satisfy equations
$$
\dot u_i^2=P_i(u_i),
$$
where
$$
P_i(u)=-(J_1+J_n)^2u^3+u^2B_i+uC_i+D_i,\qquad i=1,2;
$$
and
$$
\begin{aligned}
B_i&=2f_{i2}(J_1+J_n)-f_{i1}^2(J_1+J_n)^2,\\
C_i&=4\gamma_{(i)3}^2\chi_{(i)3}^2-f_{i2}^2+(J_1+J_n)^2f_{i1}(f_{i1}f_{i2}-2\epsilon\chi_{(i)3}\gamma_{(i)3}),\\
D_i&=-(2\epsilon \gamma_{(i)3}\chi_{(i)3}-f_{i1}f_{i2})^2,\qquad i=1,2.
\end{aligned}
$$

So,  the integration of the system
$$
\int\frac{du_1}{\sqrt{P_1(u_1)}}=t,\quad
\int\frac{du_2}{\sqrt{P_2(u_2)}}=t
$$
leads to the functions associated with the elliptic curves $E_1,
E_2$ given with:
\begin{equation}
E_i=E_i(J_1, J_n, f_{i1}, \gamma_{(i)3}, \epsilon, \chi_{(i)3},
f_{i3}): \qquad  y^2 = P_i(u). \label{8.3}
\end{equation}

The equations \eqref{8.2} are very similar to those for the
symmetric top, and they are special case of the equations of the
Lagrange bitop (see \cite{DrGa1}). From the equations \eqref{8.2}
one concludes that the dynamics of the magnetic pendulum on
$Gr^+(4,2)$ splits on two independent systems on the sphere $S^2$.
This splitting  corresponds to the fact that $Gr^+(4,2)$ is a
product of two spheres. Let us mention that a general Lagrange bitop
is more complex since it doesn't split on two independent Lagrange
tops.

\subsection{Algebro-geometric integration procedure of the magnetic pendulum}

Algebro-geometric integration completely follows paper \cite{DrGa1} (see also \cite{DrGa2}).

The starting point in the integration is the following Lax
representation:

\begin{prop}[\cite{DrGa1}]{\label{p8.1}}
The equations \eqref{8.1} has the Lax representation:
$$
\dot{L}(\lambda)=[L(\lambda), A(\lambda)]
$$
where $L(\lambda)=\lambda^2c-\lambda m+x$,
$A(\lambda)=\lambda\gamma-(J_1+J_n)m$ and
$c=\frac{1}{J_1+J_n}\gamma$
\end{prop}

\begin{rem}
The magnetic spherical pendulum on adjoint orbits \eqref{xx}
admits a similar Lax representation that provides a complete set
of commuting integrals for a regular $c\in\g$ \cite{BJ5}.
\end{rem}

We will change the coordinates in order to
diagonalize the matrix $\frac{1}{J_1+J_n}\gamma$. In this new basis the matrices
$L(\lambda)$ have the form $\tilde L (\lambda )=U^{-1} L(\lambda )U,$
$$
\tilde L(\lambda )=\left(
\begin{matrix}
-i\Delta _{34} & 0 & -\beta _3^{*} -i\beta _4^{*} & i\beta _3-\beta _4 \\
0 &  i\Delta _{34} &  -i\beta _3^{*} -\beta _4^{*} & -\beta _3+i\beta _4\\
\beta _3-i\beta _4 & -i\beta _3 +\beta _4 & -i\Delta _{12} & 0 \\
i\beta _3^{*} +\beta _4^{*} & \beta _3^{*}+i\beta _4^{*} & 0 & i\Delta _{12}
\end{matrix}\right)
$$
where
$\Delta _{12}=\lambda ^2 c_{12}-\lambda m_{12}+x _{12},\quad
\Delta _{34}=\lambda ^2c_{34}-\lambda m_{34}+x _{34},$ and
\begin{equation}
\begin{aligned}
\beta _3&=x_3+\lambda
y_3, & x_3&=\frac 12 \left( x_{13}+ix_{23}\right),\\
 \beta _4 &=
x_4+\lambda y_4, & x_4&=\frac 12 \left( x_{14}+ix_{24}\right),\\
 \beta _3^{*}&=\bar x_3+\lambda \bar y_3, & y_3&=-\frac 12 \left( m_{13}+
im_{23}\right),\\
\beta _4^{*}&=\bar x_4+\lambda \bar y_4, & y_4&=-\frac 12 \left( m_{14}+
im_{24}\right).
\end{aligned}
\label{ag2}
\end{equation}
The spectral polynomial $p(\lambda, \mu )=\det \left ( \tilde
L(\lambda )-\mu \cdot 1\right)$ has the form
\begin{equation}
p(\lambda , \mu )=\mu ^4+P(\lambda )\mu ^2 +[Q(\lambda )]^2,
\label{ag3}
\end{equation}
where
\begin{equation}
P(\lambda)=\Delta _{12}^2+\Delta _{34}^2+4\beta _3\beta _3
^{*}+4\beta _4\beta _4^{*},\quad
Q(\lambda)=\Delta _{12}\Delta _{34}+2i(\beta _3^{*}\beta _4 -
\beta _3\beta _4^{*}).
\label{ag4}
\end{equation}
We can rewrite it in terms of  $m_{ij}$ and $x_{ij}$:
\begin{equation}
P(\lambda )=A\lambda ^4 -B\lambda ^3+D\lambda ^2-E\lambda +F,\quad
Q(\lambda )=G\lambda ^4-H\lambda ^3+I\lambda ^2-J\lambda +K.
\label{ag4}
\end{equation}
Their coefficients
$$
\begin{aligned}
A&=c_{12}^2+c_{34}^2,\\
B&=2c_{34}m_{34}+2c_{12}m_{12},\\
D&=m_{13}^2+m_{14}^2+m_{23}^2+m_{12}^2+m_{34}^2+2c_{12}x_{12}+2c_{34}
x_{34},\\
E&=2x_{12}m_{12}+2x_{13}m_{13}+2x_{14}m_{14}+2x_{23}
m_{23}+2x_{24}m_{24}+2x_{34}m_{34},\\
F&=x_{12}^2+x_{13}^2+x_{14}^2+x_{23}^2+x_{24}^2
+x_{34}^2,\\
G&=c_{12}c_{34},\\
H&=c_{34}m_{12}+c_{12}m_{34},\\
I&=c_{34}x_{12}+x_{34}c_{12}+m_{12}m_{34}+m_{23}m_{14}-m_{13}m_{24}\\
J&=m_{34}x_{12}+m_{12}x_{34}+m_{14}x_{23}+m_{23}x_{14}-
x_{13}m_{24}-x_{24}m_{13},\\
K&=x_{34}x_{12}+x_{23}x_{14}-x_{13}x_{24}
\end{aligned}
$$
are integrals of the motion. From the constraints \eqref{t*o} one
can calculate values of four integrals
$$
E=2\epsilon F=2\epsilon(\chi_{12}^2+\chi_{34}^2), \qquad
J=2\epsilon K=2\epsilon\chi_{12}\chi_{34}.
$$

There is an involution
$\sigma:\;(\lambda,\mu)\rightarrow (\lambda, -\mu)$
on the curve $\Gamma: p(\lambda , \mu )=0$, which corresponds to the
skew symmetry of the matrix $L(\lambda)$. Denote the factor-curve by  $\Gamma_1=\Gamma/\sigma$.

Detailed analysis of algebro-geometric properties of the curves
$\Gamma, \Gamma_1$ one may find in \cite{DrGa1}.

We consider the next eigen-problem
$$
\left(\frac {\partial }{\partial t}+\tilde A(\lambda )\right) \psi _k=0,\quad
\tilde L(\lambda )\psi _k=\mu _k \psi _k,
$$
where $\psi _k$ are the eigenvectors with the eigenvalue
$\mu _k$. Then $\psi _k(t,\lambda )$ form $4\times 4$ matrix
with components   $\psi _k^i(t,\lambda )$. Denote by  $\varphi _i^k$
corresponding inverse matrix.

Let us introduce
$$
g_j^i(t,(\lambda , \mu _k))=\psi _k^i(t,\lambda )\cdot \varphi _j^k(t,\lambda )
$$
(there is no summation on $k$) or, in other words $g(t)=\psi
_k(t)\otimes \varphi(t) ^k.$ Matrix $g$ is of rank 1, and we have
$$
\frac{\partial \psi}{\partial t}=-\tilde{A}\psi,\qquad
\frac{\partial \varphi}{
\partial t} =\varphi \tilde{A},\qquad \frac{\partial g}{\partial t}=
[g,\tilde{A}].
$$
We can consider vector-functions $\psi
_k(t,\lambda )= \left(\psi^1_k(t, \lambda),...,\psi^4_k(t,
\lambda)\right)^{T}$ as one vector-function $\psi(t, (\lambda,
\mu))=\left(\psi^1(t,(\lambda, \mu)),..., \psi^4(t, (\lambda,
\mu))\right)^{T}$ on the  curve  $\Gamma $ defined with $\psi
^i(t,(\lambda , \mu_k))= \psi^i_k(t, \lambda)$. The same we have
for the matrix $\varphi ^k_i$. The relations for the divisors of
zeroes and poles of the functions $\psi ^i$ i $\varphi _i$ in the
affine  part of the curve $\Gamma $ are:
\begin{equation}
\left( g^i_j\right)_a=d_j(t)+d^i(t)-D_r-D'_s,
\label{gij}
\end{equation}
where  $D_r$ is the ramification divisor over $\lambda$ plane(see
\cite{Du} and \cite{DrGa1}) and $D_s$ is divisor of singular
points, $D'_s\le D_s$. One can easily calculate $\deg D_r=16,\,
\deg D_s=8$.

The matrix elements $g_j^i(t,(\lambda ,\mu _k))$  are meromorphic
functions on the curve $\Gamma $. We need their asymptotics in the
neighborhoods of the points $P_k$, which cover the point $\lambda
=\infty $.

It was justified in \cite{DrGa1} that from now on we may consider
all the functions in this section as functions on the normalization
$\tilde \Gamma $ of the curve $\Gamma $.

Let us denote by $\tilde d_j$ and by $\tilde d^i$ the following divisors:
$$
\begin{aligned}
&\tilde d_1=d_1+P_2,\quad \tilde d_2=d_2+P_1,\quad \tilde d_3=d_3+P_4,\quad
\tilde d_4 =d_4+P_3,\\
&\tilde d^1=d^1+P_2,\quad \tilde d^2=d^2+P_1,\quad \tilde d^3=d^3+P_4,\quad
\tilde d^4 =d^4+P_3.
\end{aligned}
$$

We have (see \cite{DrGa1} for details):

\begin{prop}[\cite{DrGa1}] {\rm (i)} The divisors of matrix elements of $g$ are
$$
\left( g_j^i\right)=\tilde d^i+\tilde d_j-D_r+2\left(
P_1+P_2+P_3+P_4\right)-P_i-P_j
$$

{\rm (ii)} The divisors $ \tilde d_i, \tilde d^j$ are of the same
degree
$$
\deg \tilde d_i=\deg \, \tilde d^j=5.
$$
\end{prop}

Let us denote with $\Phi(t,\lambda)$ the fundamental solution of
$$
\left(\frac{\partial}{\partial t}+\tilde{A}(\lambda)\right)\Phi(t,\lambda)=0,
$$
normalized with $\Phi(\tau)=1$.
Then, if we introduce functions
$$
\hat{\psi}^i(t,\tau,(\lambda,\mu_k))=\sum_s\Phi^i_s(t,\lambda)h^s(\tau,(\lambda,\mu_k))
$$
where $h^s$ are the eigenvector of $L(\lambda)$ normalized by the condition
$\sum_s h^s(t,(\lambda,\mu_k))=1$, it follows that
$$
\hat{\psi}^i(t,\tau,(\lambda,\mu_k))=\sum_s\Phi^i_s(t,\lambda)
\frac{\psi^s_k(\tau,\lambda)}{\sum_l\psi^l_k(\tau,\lambda)}=
\frac{\psi^i_k(t,\lambda)}{\sum_l\psi^l_k(\tau,\lambda)}.
$$

\begin{prop}[\cite{DrGa1}]{\label{p8.2}} The functions $\hat{\psi}^i$ satisfy the following
properties

{\rm (i)} In the affine part of $\tilde \Gamma$ the function
$\hat{\psi}^i$ has 4 time dependent zeroes which belong to the
divisor $d^i(t)$ defined by formula \eqref{gij}, and 8 time
independent poles, e.q.
$$
\left(\hat{\psi}^i(t,\tau,(\lambda,\mu_k))\right)_a=d^i(t)-\bar{\mathcal{D}},
\qquad \deg\bar{\mathcal{D}}=8.
$$

{\rm (ii)}  At the points $P_k$, the functions $\hat{\psi}^i$ have
essential singularities as follows:
$$
\hat{\psi}^i(t, \tau , (\lambda , \mu )) = \exp \, \left[ -(t-\tau )
R_k\right] \hat{\alpha}^i(t,\tau,(\lambda,\mu))
$$
where $R_k$ are given with
$$
R_1=i\left(\frac{\gamma_{34}}{z}-(J_1+J_n)m_{34}\right),
R_2=-R_1,
R_3=i\left(\frac{\gamma_{12}}{z}-(J_1+J_n)m_{12}\right),
R_4=-R_3
$$
and $\hat{\alpha}^i$ are holomorphic in a neighborhood of $P_k$,
$$
\hat{\alpha}^i(\tau,\tau,(\lambda,\mu))=h^i(\tau,(\lambda,\mu)),\quad
\hat{\alpha}^i(t,\tau,P_k)=\delta_i^k+v^i_k(t)z+O(z^2),
$$
with
\begin{equation}
v^i_k=\frac{\tilde{m_{ki}}}{\tilde{c_{ii}}-\tilde{c_{kk}}}.
\label{vij2}
\end{equation}
\end{prop}

We have

\begin{lem}[\cite{DrGa1}] The following relation takes place on the Jacobian
$Jac(\tilde \Gamma)$:
$$
\mathcal{A}(d^j(t)+\sigma d^j(t))=\mathcal{A}(d^j(\tau)+\sigma d^j(\tau))
$$
where $\mathcal{A}$ is the Abel map from the curve $\tilde \Gamma$ to
$Jac(\tilde \Gamma)$.
\end{lem}

From the previous Lemma we see that the vectors
$\mathcal{A}(d^i(t))$ belong to some translation of the Prym
variety $\Pi = Prym(\tilde \Gamma\,|\,\Gamma_1)$. More details
concerning the Prym varieties one can find in \cite{Mum2, Mum1,
Fay, Shok1, Shok2, DrGa1}.

It was  shown in \cite{DrGa1} that the Baker--Akhiezer function
$\Psi$ satisfies usual conditions of normalized ($n$=)4-point
function on the curve of genus $g=5$ with the divisor $\bar{
\mathcal{D}}$ of degree $\deg\bar {\mathcal{D}}= g+n-1=8$, see
\cite{DKN, Du1}. By the general theory, it should determine all
dynamics uniquely.

Let us consider the differentials
$\Omega^i_j=g_{ij}d\lambda,\quad  i,j=1,\dots , 4.$

It was proven by Dubrovin in the case of general position, that
$\Omega^i_j$ is a meromorphic differential having poles at $P_i$ and $P_j$,
with residues $v^i_j$ and $-v^j_i$ respectively. But here we have

\begin{prop} [\cite{DrGa1}] The four differentials
 $\Omega^1_2,\; \Omega^2_1,\; \Omega^3_4,\; \Omega^4_3$
 are holomorphic during the whole evolution.
 \end{prop}

\begin{rem}\label{isohol}
The proof was based on the fact that
\begin{equation}
v^1_2=v^2_1=v^3_4=v^4_3=0, \label{vij1}
\end{equation}
which is consequence of condition $\tilde L_{12}=\tilde
L_{21}=\tilde L_{34}=\tilde L_{43}=0$.  It was the reason that the
notion of {\it isoholomorphic systems} has been introduced in
\cite{DrGa1} to describe such class of integrable systems.
\end{rem}

Let us recall (\cite{DrGa1}) the general formulae for $v$:
\begin{equation}
v^i_j=\frac{\lambda_i \theta (A(P_i)-A(P_j)+tU+z_0)} {\lambda_j
\theta (tU+z_0) \epsilon (P_i,P_j)},\quad  i\ne j,
\label{vij}
\end{equation}
where $U =\sum x^{(k)}U^{(k)}$ is certain linear combination of $b$ periods $U^{(i)}$ of the differentials of the second kind
$\Omega^{(1)}_{P_i}$, which have pole of order two at $P_i$; $\lambda_i$ are
nonzero scalars, and
$$
\epsilon (P_i,P_j)=\frac {\theta [\nu ](A(P_i-P_j))}{(-\partial _{U^{(i)}}\theta [\nu ](0))^{1/2} (-\partial _{U^{(j)}}\theta [\nu ](0))^{1/2})}.
$$
(Here $\nu $ is an arbitrary odd non-degenerate characteristics.)
Thus, from \eqref{vij}, it follows

\begin{prop}[\cite{DrGa1}]  Holomorphicity of some of the differentials
$\Omega^i_j$ implies that the theta divisor of the spectral curve
contains some tori.
\end{prop}

In a case of spectral curve which is a double unramified covering
$$
\pi: \tilde \Gamma \rightarrow \Gamma _1;
$$
with $g(\Gamma _1)=g, \quad g(\tilde \Gamma )=2g-1$, as it is
satisfied for the Lagrange bitop, it is really satisfied that the
theta divisor contains a torus, see \cite{Mum2}. Following
\cite{Mum2} and \cite{DrGa1}, let us denote by $\Pi^-$ the set
$$
   \Pi^-=\left\{ L\in Pic^{2g-2} \tilde \Gamma \,|\, Nm L = K_{\Gamma 1}, h^0(L) ~
   \text {is odd} \right\},
$$
where $K_{\Gamma _1}$ is the canonical class of the curve $\Gamma
_1$ and $Nm: Pic \tilde \Gamma \rightarrow Pic \Gamma _1$ is the
norm map, see \cite{Mum2, Shok2} for details. For us,  it is crucial that
$\Pi^-$ is a translate of the Prym variety $\Pi$ and that
Mumford's relation (\cite{Mum2}, p.241-242) holds
\begin{equation}
\Pi^-\subset \Theta _{\tilde \Gamma}.
\label{mum1}
\end{equation}

Let us denote
\begin{equation}
U= i(\chi_{34} U^{(1)}-\chi_{34} U^{(2)}+ \chi_{12} U^{(3)}-\chi_{12} U^{(4)}),
\label{u1}
\end{equation}
where $U^{(i)}$ is the vector of $\tilde b$ periods of the
differential of the second kind $\Omega^{(1)}_{P_i}$, which is
normalized by the condition that $\tilde a$ periods are zero. We
suppose here that the cycles $\tilde a, \tilde b$ on the curve
$\tilde \Gamma$ and $a, b$ on $\Gamma_1$ are chosen  to correspond
to the involution $\sigma $ and the projection $\pi$, see \cite{Shok2}:
 $$
\pi (\tilde a_0)=a_0; \quad \pi(\tilde b_0)= 2 b_0,\quad
\sigma (\tilde a_k)= \tilde a_{k+2}, \quad k=1,2.
$$
The basis of normalized holomorphic differentials $[u_0,\dots ,u_5]$ on $\tilde \Gamma $ and $[v_0, v_1,v_2]$ on $\Gamma_1$  are chosen such that
$$
\pi^*(v_0)=u_0,\quad
\pi^*(v_i)=v_i+\sigma(v_i)=v_i+v_{i+2},\quad i=1,2.
$$

Now we have

\begin{thm} [\cite{DrGa1}] {\rm (i)} If the vector $z_0$ in \eqref{vij} corresponds to the translation of
the Prym variety $\Pi$ to $\Pi^-$, and the vector $U$ is defined by \eqref{u1} than
the conditions \eqref{vij1} are satisfied.

{\rm (ii)} The explicit formula for $z_0$ is
\begin{equation}
z_0=\frac{1}{2}(\hat \tau_{00},\hat \tau_{01},\hat \tau_{02},\hat \tau_{01},\hat \tau_{02}),\quad
\hat \tau_{0i}=\int_{\tilde b_0} u_i, \quad i=0, 1, 2.
\label{z0}
\end{equation}
\end{thm}

The evolution on the Jacobian of the spectral curve, as we
considered $Jac (\tilde \Gamma)$   gives the possibility to
reconstruct the evolution of the Lax matrix $L(\lambda)$ only up to
the conjugation by diagonal matrices. As it was explained in
\cite{DrGa1}, for complete integration one has to pass to the
generalized Jacobian, obtained by gluing together the infinite
points. Those points are $P_1, P_2, P_3, P_4$ and corresponding
Jacobian will be denoted as $Jac(\tilde \Gamma| \left \{ P_1, P_2,
P_3, P_4\right \}).$

It can be understood as a set of classes of relative equivalence among
the divisors on $\tilde \Gamma $ of certain degree. Two divisors of the same degree  $D_1$
and $D_2$ are called {\it equivalent relative to the points} $P_1, P_2, P_3, P_4 $,
if there exists a function $f$ meromorphic on $\tilde \Gamma $ such that
$(f)=D_1-D_2$ and $f(P_1)=f(P_2)=f(P_3)=f(P_4)$.

The generalized Abel map is defined with
$$
\tilde A(P)=(A(P),\lambda_1 (P),...,\lambda_4(P)),\quad
\lambda_i(P)=\exp\int_{P_0}^P\Omega_{P_iQ_0}, i=1,...,4,
$$
and $A(P)$ is the standard Abel map.
Here $\Omega_{P_iQ_0}$ denotes the normalized differential of the third kind,
with poles at $P_i$ and at arbitrary fixed point $Q_0$.

We will use the generalized Abel theorem as it was formulated in
\cite{DrGa1}. The generalized Jacobi inverse problem can be
formulated as the question of finding, for given $z$, points
$Q_1,\dots, Q_8$ such that
$$
\aligned
&\sum_1^8A(Q_i) - \sum_2^4A(P_i)=z+K,\\
&\lambda_j=c \exp\sum_{s=1}^8\int_{P_0}^{Q_s}\Omega _{P_jQ_0} +
\kappa_j, j=1,...4,
\endaligned
$$
 where $K$ is the Riemann constant and the constants $\kappa_j$ depend on  the curve $\tilde \Gamma$,
 the points $P_1, P_2, P_3, P_4$ and the choice of local parameters around them.

We will denote by $Q_s$ the points which belong to the divisor
$\bar {\mathcal {D}}$ from the Proposition \ref{p8.2}, and by $E$
the prime form from \cite{Fay}.  Then we have

\begin{prop}[\cite{DrGa1}]  The scalars $\lambda_j$  from the formula \eqref{vij} are
given with
$$
\lambda_j=\lambda_j^0 \exp \sum_{k\ne j}ix^{(k)}\gamma_j^k,\qquad
 \lambda_j^0=c \exp\sum_{s=1}^8\int_{P_0}^{Q_s}\Omega _{P_jQ_0} + \kappa_j,
$$
where the vector $\vec x= (x^{(1)},\dots, x^{(4)})$ denotes $t(\gamma_{34},-\gamma_{34},\gamma_{12},-\gamma_{12})$ and
$$
\gamma_i^j=\frac{d}{dk_j^{-1}}ln E(P_i,P)|_{P=P_j}.
$$
($k_j^{-1}$ is a local parameter around $P_j$.)
\end{prop}

To give the formulae for the Baker-Akhiezer function, we need some notations.
Let
$$
\alpha ^j(\vec x)=\exp[i\sum \tilde \gamma _m^jx^{(m)}]\frac {\theta (z_0)}{\theta (i\sum x^{(k)}U^{(k)} +z_0)},
$$
where
$$
\tilde \gamma ^j_m=\int _{P_0}^{P_j}\Omega^{(1)}_{P_m}, \quad m\ne j,
$$
and $\tilde \gamma^m_m$ is defined by the expansion
$$
\int _{P_0}^{P}\Omega^{(1)}_{P_m}=-k_m+\tilde \gamma^m_m + O(k_m^{-1}),\quad P \rightarrow P_m.
$$
Denote
$$
\phi ^j(\vec x, P)= \alpha ^j(\vec x) \exp(-i\int _{P_0}^P\sum x^{(m)}\Omega^{(1)}_{P_m})\frac {\theta (A(P)-A(P_j)-i\sum x^{(k)}U^{(k)}-z_0)}{\theta (A(P)-A(P_j)-z_0)}.
$$
Finally, one can state
\begin{prop} [\cite{DrGa1}] The Baker-Akhiezer function is given by
$$
\psi ^j(\vec x, P)=\phi ^j(\vec x, P)\frac {\lambda_j^0\frac {\theta (A(P-P_j)-z_0)}{\epsilon (P, P_j)}}{\sum_{k=1}^4 \lambda _k^0 \frac {\theta (A(P-P_k)-z_0)}{\epsilon (P, P_k)}}, \quad j=1,\dots ,4,
$$
where $z_0$ is given by \eqref{z0}.
\end{prop}

\section{Appendix: Basic notions of the Hamiltonian
systems}\label{sec:appendix}

\subsection{Hamiltonian systems}
Let $(P,\{\cdot,\cdot\})$ be a Poisson manifold and  $\pi$ be the
associated bivector field on $P$
$$
\{f_1,f_2\}(x)=\pi_x(df_1(x),df_2(x))=\sum_{i,j} \pi^{ij}
\frac{\partial f}{\partial x_i}\frac{\partial f_2}{\partial
x_j}\,.
$$

If $\pi$ is non-degenerate, then the two-form $\omega=\sum
\omega_{ij} dx_i \wedge dx_j$ ($\omega_{ij}\pi^{jk}=\delta_i^k$)
is a symplectic form and $(P,\omega)$ is called a symplectic
manifold.

The equations:
\begin{equation} \label{HAM}
\dot x=X_h(x) \qquad \Longleftrightarrow \qquad \dot f=\{f,h\},
\quad f\in C^\infty (P)
\end{equation}
are called {\it Hamiltonian equations} with the Hamiltonian
function $h$ and $X_h^i=\sum \pi^{ij}{\partial h}/{\partial x_j}$
is the corresponding {\it Hamiltonian vector field}.

A function $f$ is an {\it integral} of the system (constant along
trajectories of (\ref{HAM})) if and only if it commutes with $h$:
$\{h,f\}=0$. From the Jacobi identity the Poisson bracket of two
integrals is again the integral, so we can consider a Poisson
subalgebra $\mathcal F\subset C^\infty(P)$ of integrals (or a
collection of integrals closed under the Poisson bracket).
Consider the linear spaces
\begin{equation}
 F_x=\{ df(x) \, \vert\, f\in\mathcal F\}\subset T_x^* P
\label{NABLAF}
\end{equation}
and suppose that we can find  $l$ functionally independent
functions $f_1,\dots,f_l \in \mathcal F$ whose differentials span
$F_x$ almost everywhere on $M$ and that the corank of the matrix
$\{f_i,f_j\}$ is equal to some constant $k$, i.e., $\dim
\ker\pi_x\vert_{F_x}=k$. The numbers $l$ and $k$ are called {\it
differential dimension} and {\it differential index} of $\mathcal
F$ and they are denoted by $\ddim\mathcal F$ and $\dind\mathcal
F$, respectively.

We say that $\mathcal F$ is {\it complete at $x$} if the space
$F_x$ given by (\ref{NABLAF}) is coisotropic:
\begin{equation}
F^\pi_x \subset F_x \,. \label{coisotropic}
\end{equation}
Here $F^\pi_x$ is skew-orthogonal complement of $F_x$ with respect
to $\pi$:
$$
F^\pi_x=\{\xi\in T^*_x P\, \vert \,\pi_x(F_x,\xi)=0\}.
$$

The set $\mathcal F$ is {\it complete } if it is complete at a
generic point $x\in P$. In this case
$F^\pi_x=\ker\pi_x\vert_{F_x}$ and $\dind \mathcal F=\dim
F^\pi_x$, for a generic $x\in P$. Equivalently, $\mathcal F$ is
called {\it complete} if (see \cite{Bo, BJ2, Zu}):
$$
\ddim{\mathcal F}+\dind{\mathcal F}=\dim P+\corank\{\cdot,\cdot\}.
$$

The Hamiltonian system on (\ref{HAM}) is {\it completely
integrable (in noncommutative sense)} if it possesses a complete
set of first integrals $\mathcal F$. Then (under compactness
condition) $P$ is almost everywhere foliated by $(\dind{\mathcal
F}-\corank\{\cdot,\cdot\}$)-dimensional invariant tori. As in the
Liouville-Arnol'd theorem \cite{Ar}, the Hamiltonian flow
restricted to regular invariant tori is quasi-periodic (see
Nekhoroshev \cite{N}, Mishchenko and Fomenko \cite{MF2} and Zung
\cite{Zu}).

\subsection{Natural mechanical systems}
The basic examples of Hamiltonian systems are natural mechanical
systems $(Q,\kappa,v)$, where $Q$ is a configuration space,
$\kappa$ is a Riemannian metric on $Q$ and $v:Q\to\R$ is a
potential function. Let $q=(q^1,\dots,q^n)$ be local coordinates
on $Q$. The motion of the system is described by the
Euler--Lagrange equations
\begin{equation} \label{Lagrange}
\frac{d}{dt}\frac{\partial l}{\partial \dot q^i}=\frac{\partial
l}{\partial q^i}, \quad i=1,\dots,n,
\end{equation}
where the Lagrangian is $l(q,\dot q)=\frac12(\kappa_q \dot q,\dot
q)-v(q) =\frac12\sum_{ij}\kappa_{ij}\dot q^i\dot q^j-v(q)$.

Equivalently, we can pass from velocities $\dot q^i$ to the
momenta $p_j$ by using the standard Legendre transformation
$p_j=\kappa_{ij}\dot q^i$. Then in the coordinates $q^i, p_i$ of
the cotangent bundle $T^*Q$ the equations of motion read:
\begin{equation} \label{1}
\frac{dq^i}{dt}=\frac{\partial h}{\partial p_i},\qquad
\frac{dp_i}{dt}=-\frac{\partial h}{\partial q^i}, \qquad
i=1,\dots,n,
\end{equation}
where the Hamiltonian $h$ is the sum of the kinetic and potential
energy of the system $
 h(q,p)=\frac{1}{2}\sum_{i,j}\kappa^{ij}p_i p_j+v(q).
$ Here $\kappa^{ij}$ are the coefficients of the tensor inverse to
the metric.

This system of equations is Hamiltonian on $T^*Q$ endowed with the
\emph{canonical symplectic form} $\omega=\sum_{i=1}^n dp_i\wedge
dq^i$. The corresponding canonical Poisson bracket is given by
\begin{equation}
\{f,g\}=\sum_{i=1}^{n}\left(
 \frac{\partial f}{\partial q^i}\frac{\partial g}{\partial p_i}
-\frac{\partial g}{\partial q^i}\frac{\partial f}{\partial
p_i}\right). \label{CPB}\end{equation}

Let $\epsilon$ be a real parameter (a "coupling" constant). The
motion of the particle under the influence of the additional
magnetic field given by a closed 2-form $ \epsilon\,\Omega=\sum_{1
\le i<j \le n} \epsilon \,F_{ij}(q) dq^i \wedge dq^j, $ is
described by the following equations:
\begin{equation}
\frac{dq^i}{dt}=\frac{\partial h}{\partial p_i}, \qquad
\frac{dp_i}{dt}=-\frac{\partial h}{\partial
q^i}+\epsilon\sum_{j=1}^{n} F_{ij} \frac{\partial H}{\partial
p_j}. \label{magnetic_flow}
\end{equation}

The equations (\ref{magnetic_flow}) are Hamiltonian with respect
to the symplectic form $\omega+\epsilon \rho^*\Omega$, where
$\rho: T^*Q\to Q$ is the natural projection. Namely, the new
Poisson bracket is given by
\begin{equation}
\{f,g\}_\epsilon=\{f,g\}+\epsilon\sum_{i,j=1}^{n} F_{ij}
\frac{\partial f}{\partial p_i}\frac{\partial g}{\partial p_j},
\label{magnetic_bracket}
\end{equation}
and the Hamiltonian equations $\dot f=\{f,h\}_\epsilon$ read
(\ref{magnetic_flow}).

\subsection{Hamiltonian $G$-actions}
Let a connected Lie group $G$ act on $2n$-dimen\-sional connected
symplectic manifold $(M,\omega)$. The action is \emph{Hamiltonian}
if $G$ acts on $M$ by symplectomorphisms and there is a
well-defined momentum mapping:
\begin{equation} \label{moment_map}
\Phi: M \to \mathfrak g^*
\end{equation}
($\mathfrak g^*$ is a dual space of the Lie algebra $\mathfrak g$)
such that one-parameter subgroups of symplectomorphisms are
generated by the Hamiltonian vector fields of functions
$\phi_\xi(y)=(\Phi(y),\xi)$, $\xi\in\mathfrak g$ and $
\phi_{[\xi_1,\xi_2]}=\{\phi_{\xi_1},\phi_{\xi_2}\}. $ Then $\Phi$
is equivariant with respect to the given action of $G$ on $M$ and
the co-adjoint action of $G$ on $\mathfrak g^*$: $ \Phi(g\cdot
x)=\Ad_g^*(\Phi(x)). $ In particular, if $\eta$ belongs to
$\Phi(M)$, then the co-adjoint orbit ${\mathcal O}(\eta)$ belongs
to $\Phi(M)$ as well.

The mapping $f\mapsto f\circ\Phi$ is a morphism of Poisson
structures: $
\{f_1\circ\Phi,f_2\circ\Phi\}(x)=\{f_1,f_2\}_{\mathfrak
g^*}(\eta)$,  $\eta=\Phi(x)$, where $\{\cdot,\cdot\}_{\mathfrak
g^*}$ is the Lie--Poisson bracket on $\mathfrak g^*$:
\begin{equation*} \label{2.5}
\{f_1,f_2\}_{\mathfrak g^*}(\eta)=(\eta,[df_1(\eta),df_2(\eta)]),
\quad f_1,f_2: \mathfrak g^*\to \mathbb{R}.
\end{equation*}

Thus, $\Phi^*C^\infty(\mathfrak g^*)$ is closed under the Poisson
bracket. Since $G$ acts in a Hamiltonian way, the set of
$G$-invariant functions $C^\infty_G(M)$ in $C^\infty(M)$ is closed
under the Poisson bracket as well. Also
$\{\Phi^*C^\infty(\mathfrak k^*),C^\infty_G(M)\}=0$ (the Noether
theorem).

Suppose the group $G$ is compact. Then we have

\begin{thm} [\cite{BJ2}] \label{Bol_Jov}
\emph{(i)} The algebra of functions $\Phi^*C^\infty(\mathfrak
g^*)+C^\infty_G(M)$ is complete:
$$
\ddim(\Phi^*C^\infty(\mathfrak
g^*)+C^\infty_G(M))+\dind(\Psi^*C^\infty(\mathfrak
g^*)+C^\infty_G(M))=\dim M.
$$

\emph{(ii)} Suppose $\mathcal A\subset C^\infty (\mathfrak g^*)$
is a involutive set of functions, complete on a generic coadjoint
orbit ${\mathcal O}(\eta)\subset \Phi(M)$ and $\mathcal B$ is a
complete commutative subset of $C^\infty_G(M)$:
$$
\ddim\mathcal B=\dind\mathcal B=\frac12\left(\ddim
C^\infty_G(M)+\dind C^\infty_G(M)\right).
$$
Then $\Phi^*{\mathcal A}+C^\infty_G(M)$ and
$\Phi^*C^\infty(\mathfrak g^*)+\mathcal B$ are complete sets on
$M$, while $\Phi^*{\mathcal A}+{\mathcal B}$ is a complete
commutative set on $M$.
\end{thm}

\subsection{Symplectic reductions}\label{sim-red}
Let $G$ be a Lie group with a free and proper Hamiltonian action
on a symplectic manifold $(M,\omega)$ with the momentum mapping
(\ref{moment_map}). Assume that $\eta$ is a regular value of
$\Phi$, so that $M_\eta=\Phi^{-1}(\eta)$ and $M_{\mathcal
O_\eta}=\Phi^{-1}(\mathcal O_\eta)$ are smooth manifolds. Here
$\mathcal O_\eta=G/G_\eta$ is the coadjoint orbit of $\eta$. The
manifolds $M_\eta$ and $M_{\mathcal O_\eta}$ are
$G_\eta$-invariant and $G$-invariant, respectively. There is a
unique symplectic structure $\omega_\eta$ on $N_\eta=M_\eta/G_\eta
\cong \Phi^{-1}(\mathcal O_\eta)/G$ satisfying
$$
\omega|_{M_\eta}=d\pi_\eta^*\omega_\eta,
$$
where $\pi_\eta: M_\eta\to N_\eta$ is the natural projection
(Marsden and Weinstein \cite{MaWe}). According to Noether's
theorem, if $h$ is a $G$-invariant function, then the momentum
mapping $\Phi$ is an integral of the Hamiltonian system $\dot
x=X_h(x)$. In addition, the restriction of $X_h$ to the invariant
submanifold $M_\eta$ projects to the Hamiltonian vector field
$X_{h_\eta}$ on the reduced space $N_\eta$ with $h_\eta$ defined
by $h|_{M_\eta}=\pi_\eta^* h_\eta=h_\eta\circ\pi_\eta$.

\subsection{Cotangent bundle reductions}
As the important example, consider the $G$--action on the
configuration space $Q$. The action can be naturally extended to
the Hamiltonian action on $(T^*Q,\omega)$: $ g\cdot (q,p)=(g\cdot
q, (dg^{-1})^*p) $ with the momentum mapping $\Phi$ given by
\begin{equation} \label{ctg-mm}
(\Phi(q,p),\, \xi)=(p,\, \xi_q), \quad\xi\in\g,
\end{equation}
where $\xi_q$ is the vector given by the action of one-parameter
subgroup $\exp(t\xi)$ at $q$ \cite{LM}.

Now, let $G$ be a connected Lie group acting freely and properly
on $Q$ and $\pi: Q\to Q/G$ be the canonical projection. Then $0$
is the regular value of the cotangent bundle momentum mapping
(\ref{ctg-mm}) and the reduced space $(\Phi^{-1}(0)/G,\omega_0)$
is symplectomorphic to $T^*(Q/G)$.

Suppose $(Q,\kappa,v)$ is a $G$-invariant natural mechanical
system. That is $G$ acts by isometries and the potential is the
pull back of the potential $V$ defined on $Q/G$. The metric
$\kappa$ induce the submersion metric on $Q/G$ (e.g., see
\cite{Be}). The reduced system, for a zero value of the momentum
mapping, is the natural mechanical system $(Q/G,K,V)$. For Abelian
groups this is the classical method of Routh for eliminating
cyclic coordinates \cite{Ro}. Within Lagrangian formalism the
non-Abelian construction for the zero level-set of and for the
other values of the momentum mapping are given in \cite{AKN} and
\cite{MRS}, respectively.

\subsection{Compatible Poisson brackets}
Let $\{\cdot,\cdot\}_1$ and $\{\cdot,\cdot\}_2$ be compatible
Poisson structures on a manifold $P$. In other words, each linear
combination  $\lambda_1 \{\cdot,\cdot\}_1+ \lambda_2
\{\cdot,\cdot\}_2$ with constant coefficients is again a Poisson
structure (e.g., see \cite{Bo, TF, Za, Pa} and references there
in). Let $\pi_1$ and $\pi_2$ be the associated bivector fields and
let
$$ \Pi=\{\pi_{\lambda_1,\lambda_2} \; \vert\;
\lambda_1,\lambda_2\in \mathbb{R}, \;\lambda_1^2+\lambda_2^2\ne
0\}, \quad \pi_{\lambda_1,\lambda_2}=\lambda_1 \pi_1+\lambda_2
\pi_2.
$$

In what follows we shall suppose that all functions are defined on
some open set $U$, $x\in U$. By $r$ denote the corank of a generic
bracket (or bivector)  in $\Pi$ at $x$. For each bracket in $\Pi$
of corank $r$, we consider the set of its Casimir functions at
$x$. Let ${\mathcal C}$ be the union of these sets. Then
${\mathcal C}$ is involutive set with respect to every Poisson
bracket from $\Pi$. Let $C_x$ denote the linear subspace of $T^*_x
P$ generated by the differentials of functions from ${\mathcal
C}$. It is clear that $C_x$ is spanned by the kernels $\ker
\pi(x)$, $\pi\in\Pi$, $\corank\pi(x)=r$.

Together with $\Pi$, consider its natural complexification
$\Pi^{\mathbb{C}}=\{\pi_{\lambda_1,\lambda_2}=\pi_1\Lambda_1+
\pi_2\Lambda_2$, $\lambda_1,\lambda_2\in \mathbb{C}, \;
\vert\lambda_1\vert^2+\vert\lambda_2\vert^2\ne 0\}$. Here, we
consider $\pi_{\lambda_1,\lambda_2}$ as a complex valued
skew-symmetric bilinear form on the complexification of the
co-tangent space $(T^*_x P)^\mathbb{C}$. There are only finite
number of the non-proportional singular structures
$\pi_{\lambda_1^1,\lambda_2^1},\dots,\pi_{\lambda_1^\rho,\lambda_2^\rho}\in
\Pi^\mathbb{C}$ with a corank greater then $r$ at $x$.
 With the above
notation, we can state the following remarkable result:

\begin{thm}[Bolsinov~\cite{Bo}] \label{bolsinov} {\rm(i)} $C^\pi_x$ does not depend on
the choice $\pi\in\Pi$.

{{\rm (ii)}} $(C^\pi_x)^\mathbb{C}\supset
C_x^\mathbb{C}+\ker\pi_{\lambda_1^1,\lambda_2^1}(x)+\dots+\ker\pi_{\lambda_1^\rho,\lambda_2^\rho(x)}$.

{{\rm (iii)}} $(C^\pi_x)^\mathbb{C}=C_x^\mathbb{C}+
\ker\pi_{\lambda_1^1,\lambda_2^1}(x)+\dots+\ker\pi_{\lambda_1^\rho,\lambda_2^\rho(x)}$
if and only if
\begin{equation}
\dim_\mathbb{C} K_{\lambda_1^i,\lambda_2^i}=r, \quad
K_{\lambda_1^i,\lambda_2^i}=\ker \pi_0\vert_{\ker
\pi_{\lambda_1^i,\lambda_2^i}}\subset (T_x^*C)^\mathbb{C}, \quad
i=1,\dots,\rho,\label{uslovi-bolsinova}
\end{equation}
where $\pi_0\in \Pi$ is of the maximal rank at $x$.
\end{thm}

As a corollary, an important completeness condition is formulated
in \cite{Bo}:

\begin{thm} [\cite{Bo}]\label{bolsinov*}
Let $\pi\in \Pi$ and $\corank \pi(x)=r$. Then ${\mathcal C}$ is a
complete commutative set at $x \in P$ if and only if $\corank
\pi'(x)=r$ for all $\pi'\in\Pi^\mathbb{C}$, $\pi'\ne\lambda\pi$,
$\lambda\in\mathbb{C}$.
\end{thm}

By the use of Theorem \ref{bolsinov} one can also formulate
conditions for non-commutative integrability in the case that some
of the brackets in $\Pi$ are not of the maximal rank \cite{Bo}
(see also \cite{Pa}).

\subsection*{Acknowledgments}
The research was supported by the Serbian Ministry of Science
Project 144014 Geometry and Topology of Manifolds and Integrable
Dynamical Systems.

\end{document}